\documentclass{amsart} 
\usepackage{amsmath}
\usepackage{amsthm}
\usepackage{hyperref}

\usepackage{amsfonts}
\usepackage{amssymb}
\usepackage{graphicx}
\usepackage{bbm}
\usepackage{tikz}
\usetikzlibrary{calc}
\usepackage[all]{xy}
\usepackage{blkarray}
\usepackage{indentfirst}
\usepackage{mathtools}
\usepackage{caption}
\usepackage{subcaption}
\newtheorem{theorem}{Theorem}[section]

\newtheorem{proposition}{Proposition}[section]

\newcommand{\Ad}{\mathrm{Ad}}
\newcommand{\ad}{\mathrm{ad}}

\newcommand{\CC}{\mathbb{C}}
\newcommand{\RR}{\mathbb{R}}

\newcommand{\Ind}{\mathrm{Ind}}

\newcommand{\ZZ}{\mathbb{Z}}

\newcommand{\Tr}{\mathrm{Tr}}

\newcommand{\frakk}{\mathfrak{k}}
\newcommand{\fraka}{\mathfrak{a}}
\newcommand{\frakg}{\mathfrak{g}}
\newcommand{\frakt}{\mathfrak{t}}
\newcommand{\frakh}{\mathfrak{h}}
\newcommand{\frakp}{\mathfrak{p}}

\newcommand{\fraku}{\mathfrak{u}}
\newcommand{\frakn}{\mathfrak{n}}
\newcommand{\frakm}{\mathfrak{m}}
\newcommand{\fraks}{\mathfrak{s}}

\newcommand{\fraksu}{\mathfrak{su}}
\newcommand{\fraksl}{\mathfrak{sl}}

\newcommand{\Sym}{\mathrm{Sym}}
\newcommand{\diag}{\mathrm{diag}}

\newcommand{\ii}{\mathsf{i}}
\newcommand{\dr}{\mathsf{d}r}
\newcommand{\dd}{\mathsf{d}}
\newcommand{\dl}{\mathsf{d}l}

\newcommand{\ppp}{\mathsf{p}}
\newcommand{\qqq}{\mathsf{q}}

\newenvironment{spmatrix}{\left(\begin{smallmatrix}}{\end{smallmatrix}\right)}

\begin{document} 

\begin{abstract}
In this paper, following a similar procedure developed by Buttcane and Miller in \cite{MillerButtcane} for $SL(3,\RR)$, the $(\frakg,K)$-module structure of the minimal principal series of real reductive Lie groups $SU(2,1)$ is described explicitly by realizing the representations in the space of $K$-finite functions on $U(2)$. Moreover, by combining combinatorial techniques and contour integrations, this paper introduces a method of calculating intertwining operators on the principal series. Upon restriction to each $K$-type, the matrix entries of intertwining operators are represented by $\Gamma$-functions and Laurent series coefficients of hypergeometric series. The calculation of the $(\frakg,K)$-module structure of principal series can be generalized to real reductive Lie groups whose maximal compact subgroup is a product of $SU(2)$'s and $U(1)$'s.\end{abstract}
 
\title[$SU(2,1)$ Principal Series]{Principal Series Representation of $SU(2,1)$ and Its Intertwining Operator} 
\author{Zhuohui Zhang}
\email{zhuohui.zhang@weizmann.ac.il}
%\thanks{Stephen D Miller}
%\affiliation{Faculty of Mathematics and Computer Science, Weizmann Institute of Science, Israel}
\maketitle

\section{Introduction} 
\subsection{Automorphic Forms and $(\frakg,K)$-Modules}
The study of automorphic forms serves as a central topic in representation theory and number theory. The \emph{modular forms} and \emph{Maa{\ss} forms} on the upper half plane
\[
	\mathbb{H} = \{x+\ii y|x,y\in\RR, y>0\}\subset \CC
\]
under the action of an arithmetic subgroup $\Gamma$ of $SL(2,\ZZ)$ acting on $\mathbb{H}$ by fractional linear transform $\begin{spmatrix}a&b\\c&d\end{spmatrix}\cdot z\mapsto \frac{az+b}{cz+d}$
are two classical objects in this area, connecting the study of algebraic curves, representation theory and number theory. They are defined as certain eigenfunctions of the weight $k\in \ZZ^{\geq 0}$ Laplace operator $\Delta_k = -y^2(\partial^2_x+\partial^2_y)+\ii ky\partial_x$, and invariant under the action of $\Gamma$ in the sense
\begin{equation}\label{slashcond}
	f|_k{\begin{spmatrix}
a&b\\c&d
\end{spmatrix}}(z):=\left(\frac{c z+d}{|c z+d|}\right)^{-k}f\left(\frac{az+b}{cz+d}\right)=f(z).
\end{equation}

The study of automorphic forms in general passes the function from the upper half plane $\mathbb{H}$ to the real reductive Lie group $G = SL(2,\RR)$. $G$ has a maximal compact subgroup $K=SO(2,\RR)=\{\begin{spmatrix}
    	\cos\theta&\sin\theta\\-\sin\theta&\cos\theta
    \end{spmatrix}\vert\theta\in\RR\}$, a real split Cartan subgroup $A =\{ \begin{spmatrix}\sqrt{y}&0\\0&1/\sqrt{y}\end{spmatrix}|y>0\}$ and a nilpotent subgroup $N = \{\begin{spmatrix}
    1&x\\0&1
    \end{spmatrix}\vert x\in\RR\}$. These subgroups give rise to an \emph{Iwasawa decomposition}
    \[\begin{spmatrix}
    	a&b\\c&d
   	\end{spmatrix}=\begin{spmatrix}
    	1&\frac{a c+b d}{c^2+d^2}\\
        0&1
    	\end{spmatrix}\begin{spmatrix}
    	1/\sqrt{c^2+d^2}&0\\
        0&\sqrt{c^2+d^2}
    	\end{spmatrix}\begin{spmatrix}
    	\frac{d}{\sqrt{c^2+d^2}}&-\frac{c}{\sqrt{c^2+d^2}}\\
        \frac{c}{\sqrt{c^2+d^2}}&\frac{d}{\sqrt{c^2+d^2}}
    	\end{spmatrix}\]
of any element of $G$. The Iwasawa decomposition $G = NAK$ parametrizes any element in $G$ with the coordinates \[(x,y,\theta) = \left(\frac{a c+b d}{c^2+d^2},1/(c^2+d^2),\arctan(-c/d)\right).\]
The upper half plane $\mathbb{H}$ is thus isomorphic to the hermitian symmetric space $G/K$, on which the fractional linear action by $N$ is the translation along $x$-axis, and the action by $A$ is the positive scalar multiple of a point.\\
\indent Thus we can consider the spaces of weight $k$ automorphic functions
\[C^{\infty}(\Gamma\backslash G,k)=\left\{f:G\rightarrow\CC\textrm{ smooth}|f(\gamma g \begin{spmatrix}
	\cos\theta &\sin\theta\\-\sin\theta &\cos\theta
	\end{spmatrix})=e^{\ii k\theta}f(g)\text{ for any }\gamma\in\Gamma\right\}\]
on which the group $G$ acts by the right regular action $\pi(g)f(x)=f(xg)$. Any weight $k$ Maa{\ss} form $f$ defines a function $F\in C^\infty(\Gamma\backslash G,k)$ via the formula $F(g) = (f|_k g)(\ii) $. Conversely, any function $F \in C^\infty(\Gamma\backslash G,k)$ defines a function $f(x+\ii y) = F\begin{spmatrix}\sqrt{y}&x/\sqrt{y}\\0&1/\sqrt{y}\end{spmatrix}$ on the upper half plane $\mathbb{H}$, which satisfies the same invariance condition under $\Gamma$ as (\ref{slashcond}). This correspondence between weight $k$ automorphic forms on $\mathbb{H}$ and the $\Gamma$-left invariant functions on $G$ on which $K$ acts on the right as a character $e^{\ii k\theta}$ motivates the study of irreducible representations of $G$ and their \emph{$K$-types}. The concept to study is the \emph{$(\frakg,K)$-module} or \emph{Harish-Chandra module} introduced by Harish-Chandra and James Lepowsky in \cite{Lepowsky}. 

\subsection{Bargmann's Classification of $SL(2,\RR)$ and $GL(2,\RR)$ Irreducible $(\frakg,K)$-Modules}\label{SL2RTheory}
We define the \emph{principal series} to be the set of $K$-finite smooth functions
%    \[
%    	I(\chi_{\delta,\lambda}) = \left\{f:G\rightarrow\CC|f\left((-1)^{\epsilon}g\begin{spmatrix}\sqrt{y}&0\\0&1/\sqrt{y}\end{spmatrix}\begin{spmatrix}1&x\\0&1\end{spmatrix}\right)=y^{-\frac{\lambda+1}{2}}(-1)^{\epsilon\delta}f(g)\right\}
%    \]
        \[
    	I(\chi_{\delta,\lambda}) = \left\{f:G\rightarrow\CC|f\left((-1)^{\epsilon}\begin{spmatrix}1&x\\0&1\end{spmatrix}\begin{spmatrix}\sqrt{y}&0\\0&1/\sqrt{y}\end{spmatrix}g\right)=y^{\frac{\lambda+1}{2}}(-1)^{\epsilon\delta}f(g)\right\}
    \]
    where $G$ acts by translation on the right. Since the group $G$ has an Iwasawa decomposition $G = NAK$, the value of $f$ is determined by its restriction to $K$. We can expand $f\in I(\chi_{\delta,\lambda})$ into finite Fourier series
    \[
    	I(\chi_{\delta,\lambda}) = \bigoplus_{k\equiv \delta\text{ }\mathrm{mod}\text{ } 2}\CC e^{\ii k \theta}.
    \]
     The irreducible $K$-representations $\CC e^{\ii k\theta}$ contained in the $(\frakg,K)$-module of $I(\chi_{\epsilon,\lambda})$ are called \emph{$K$-types}. A representation of $G$ is called \emph{admissible} if all $K$-types occur with finite multiplicities. According to Bargmann \cite{Bargmann}, the irreducible $(\frakg,K)$-modules of admissible representations of $SL(2,\RR)$ are classified by the following theorem:
\begin{theorem}\cite{Bargmann}\cite{KnappRank1}\cite{Muic}
	The $(\frakg,K)$-modules of irreducible admissible representations of $SL(2,\RR)$ can be realized as subrepresentations or quotient representations of the principal series $I(\chi_{\delta,\lambda})$ as follows:
	\begin{enumerate}
		\item If $\lambda+1\equiv \delta\text{ }\mathrm{mod}\text{ }2\ZZ$, 
			\begin{enumerate}
				\item If $\lambda>0$, $I(\chi_{\delta,\lambda})$ has two irreducible subrepresentations $D^\pm_{\lambda}$ called the discrete series representations. The quotient $W_{\lambda} = I(\chi_{\delta,\lambda})/(D^+_{\lambda}\oplus D^-_{\lambda})$ has finite dimension $\lambda$.
				\item If $\lambda<0$, $I(\chi_{\delta,\lambda})$ has a finite dimensional subrepresentation $W_{-\lambda}$ of dimension $-\lambda$. The quotient $I(\chi_{\delta,\lambda})/W_{-\lambda} \cong D^+_{-\lambda}\oplus  D^-_{-\lambda}$ splits into a direct sum of two discrete series representations.
				\item If $\delta = 1$ and $\lambda=0$, then the principal series decomposes into two limits of discrete series, and $I(\chi_{-1,0}) = D^+_{0}\oplus  D^-_{0}$.
			\end{enumerate}
		\item In all other cases,  $I(\chi_{\delta,\lambda})$ is irreducible, and $I(\chi_{\delta,\lambda})$ is isomorphic to $I(\chi_{\delta,-\lambda})$.
	\end{enumerate}
	For an arbitrary $\lambda\in\ZZ$, the (limit of) discrete series representations $D^\pm_{|\lambda|}$ have a decomposition into $K$-types
	\[
		D^\pm_{|\lambda|} = \bigoplus_{\substack{k\geq |\lambda|+1\\k\equiv\delta\text{ }\mathrm{mod}\text{ } 2}}\CC e^{\pm\ii k\theta}.
	\]
\end{theorem}
The \emph{raising and lowering operators} 
%$U_+ = \frac{1}{2}\begin{spmatrix}
%    1&-\ii\\-\ii&-1
%    \end{spmatrix}$ and $U_- = \frac{1}{2}\begin{spmatrix}
%    1&\ii\\\ii&-1
%    \end{spmatrix}$,
$U_+ = \frac{1}{2}\begin{spmatrix}
    1&\ii\\\ii&-1
    \end{spmatrix}$ and $U_- = \frac{1}{2}\begin{spmatrix}
    1&-\ii\\-\ii&-1
    \end{spmatrix}$,
      act on the $K$-types by the formula \[
    	U_\pm e^{\ii k\theta} = \frac{\lambda\pm k+1}{2}e^{\ii (k\pm 2)\theta}.
    \]
For example, the decomposition of the principal series $I(\chi_{0,1})$ can be displayed in the following diagram:
\[
\xymatrix{
         e^{-\ii k\theta}\ar@/^/[r]^{U_+} &\ldots \ar@/^/[r]^{U_+} \ar@/^/[l]^{U_-}&e^{-2\ii\theta}  \ar@/^/[l]^{U_-}&1 \ar@/^/[r]^{U_+} \ar@/^/[l]^{U_-}&e^{2\ii\theta} \ar@/^/[r]^{U_+} &\ldots \ar@/^/[r]^{U_+} \ar@/^/[l]^{U_-}&e^{\ii k\theta} \ar@/^/[l]^{U_-}
        }
\]

In the $GL(2,\RR)$ case, consider $\delta_1,\delta_2\in \{\pm 1\}$ and $\lambda_1,\lambda_2\in\CC$, and define the character $\chi_{\delta_1,\lambda_1}\times\chi_{\delta_2,\lambda_2}$ on the Cartan subgroup \[H = \left\{\begin{spmatrix}\epsilon_1&0\\0&\epsilon_2\end{spmatrix}\begin{spmatrix}a_1&0\\0&a_2\end{spmatrix}| \epsilon_1,\epsilon_2\in\{\pm 1\}, a_1,a_2\in \RR_{> 0}\right\}\] by sending the elements $\begin{spmatrix}\epsilon_1&0\\0&\epsilon_2\end{spmatrix}\begin{spmatrix}a_1&0\\0&a_2\end{spmatrix}$ to $\epsilon_1^{\delta_1}\epsilon_2^{\delta_2}a_1^{\lambda_1}a_2^{\lambda_2}$. We define the principal series $I(\chi_{\delta_1,\lambda_1}\times\chi_{\delta_2,\lambda_2})$ for $GL(2,\RR)$ to be the induced representation of $\chi_{\delta_1,\lambda_1}\times\chi_{\delta_2,\lambda_2}$ from the Borel subgroup $B = \left\{\begin{spmatrix}a&b\\0&d\end{spmatrix}|a d\neq 0\right\}$ to $G$:
%\begin{align}
%	I(\chi_{\delta_1,\lambda_1}\times\chi_{\delta_2,\lambda_2})=\{f:G\longrightarrow \CC|f(g\begin{spmatrix}\epsilon_1&0\\0&\epsilon_2\end{spmatrix}\begin{spmatrix}a_1&b\\0&a_2\end{spmatrix}) = \epsilon_1^{\delta_1}\epsilon_2^{\delta_2}a_1^{-\lambda_1-\frac{1}{2}}a_2^{-\lambda_2+\frac{1}{2}}f(g)\}
%\end{align}
\begin{align}
	I(\chi_{\delta_1,\lambda_1}\times\chi_{\delta_2,\lambda_2})=\{f:G\longrightarrow \CC|f(\begin{spmatrix}\epsilon_1&0\\0&\epsilon_2\end{spmatrix}\begin{spmatrix}a_1&b\\0&a_2\end{spmatrix}g) = \epsilon_1^{\delta_1}\epsilon_2^{\delta_2}a_1^{\lambda_1+\frac{1}{2}}a_2^{\lambda_2-\frac{1}{2}}f(g)\}
\end{align}
The description of the principal series and the classification of $GL(2,\RR)$ representations can be summarized as in Theorem 2.4 of \cite{Muic} as follows:
\begin{theorem}
	The $(\frakg,K)$-modules of irreducible admissible representations for the group $GL(2,\RR)$ can be realized as sub- or quotient modules of the principal series $I(\chi_{\delta_1,\lambda_1}\times\chi_{\delta_2,\lambda_2})$. If we define
	\[
		\begin{matrix}s=\frac{\lambda_1-\lambda_2+1}{2}, &\mu = \frac{\lambda_1+\lambda_2}{2}, &\delta = \delta_1+\delta_2\end{matrix}
	\]
	\begin{enumerate}
		\item If $s\notin \{\frac{k}{2}| k\in\ZZ, k\equiv \delta\text{ }\mathrm{mod}\text{ } 2\ZZ\}$, then $I(\chi_{\delta_1,\lambda_1}\times\chi_{\delta_2,\lambda_2})$ is irreducible. Moreover, if we interchange the two induction parameters, $I(\chi_{\delta_1,\lambda_1}\times\chi_{\delta_2,\lambda_2})\cong I(\chi_{\delta_2,\lambda_2}\times\chi_{\delta_1,\lambda_1})$.
		\item If $\lambda_1 >\lambda_2$ and $ \lambda_1-\lambda_2+1\equiv \delta\text{ }\mathrm{mod}\text{ } 2\ZZ$,  then if we set $k = \lambda_1-\lambda_2+1$, the character $\chi_{\delta_1,\lambda_1}\times\chi_{\delta_2,\lambda_2}$ takes the form $(\chi_{\delta_0,\mu}\cdot\chi_{k,\frac{k-1}{2}})\times(\chi_{\delta_0,\mu}\cdot\chi_{0,-\frac{k-1}{2}})$ where $\chi_{\delta,\lambda}$ is the character on $\RR^\times$ as defined above, sending each $a\in\RR^\times$ to $\mathrm{sgn}(a)^{\delta}|a|^\lambda$, and $\delta_0\in\{\pm 1\}$,
	then there exists a composition series for the principal series $I(\chi_{\delta_1,\lambda_1}\times\chi_{\delta_2,\lambda_2})$:
		\begin{align*}
			D^{\chi_{\delta_0,\mu}}_{k}\hookrightarrow  &I\left((\chi_{\delta_0,\mu}\cdot\chi_{k,\frac{k-1}{2}})\times(\chi_{\delta_0,\mu}\cdot\chi_{0,-\frac{k-1}{2}})\right) \rightarrow W^{\chi_{\delta_0,\mu}}_{k}
		\end{align*}
		where $D^{\chi_{\delta_0,\mu}}_{k}$ is a discrete series representations, and $W^{\chi_{\delta_0,\mu}}_{k}$ is a finite dimensional dimension. The superscript $\chi_{\delta_0,\mu}$ indicates that the center of $GL(2,\RR)$ acts by a character $\chi_{\delta_0,\mu}$.
		\item If $\lambda_1 <\lambda_2$ and $ \lambda_1-\lambda_2+1\equiv \delta\text{ }\mathrm{mod}\text{ } 2\ZZ$, then if we set $k = \lambda_2-\lambda_1+1$, the character $\chi_{\delta_1,\lambda_1}\times\chi_{\delta_2,\lambda_2}$ takes the form $(\chi_{\delta_0,\mu}\cdot\chi_{k,-\frac{k-1}{2}})\times(\chi_{\delta_0,\mu}\cdot\chi_{0,\frac{k-1}{2}})$, and the principal series representation $I\left(\chi_{\delta_1,\lambda_1}\times\chi_{\delta_2,\lambda_2}\right)$ has a composition series:
		\begin{align*}
			W^{\chi_{\delta_0,\mu}}_{k}\hookrightarrow  &I\left((\chi_{\delta_0,\mu}\cdot\chi_{k,\frac{k-1}{2}})\times(\chi_{\delta_0,\mu}\cdot\chi_{0,-\frac{k-1}{2}})\right) \rightarrow D^{\chi_{\delta_0,\mu}}_{k}.
		\end{align*}
		The notations for specific representations are the same as in $\lambda_1 > \lambda_2$ case.
	\end{enumerate}
	Denote by $\sigma_{l} = \Ind_{SO(2)}^{O(2)} \CC e^{\ii l\theta}$ the index 2 induction of a character from $SO(2)$ to $K=O(2)$. Since the element $\begin{spmatrix}1&0\\0&-1\end{spmatrix}\in O(2)$, and the conjugation of $\begin{spmatrix}\cos\theta&\sin\theta\\-\sin\theta&\cos\theta\end{spmatrix}$ by this element sends $e^{\ii l\theta}$ to $e^{-\ii l\theta}$. It is easy to see that $\sigma_l\cong \sigma_{-l}$. The $K$-type decomposition of the two irreducible $(\frakg,K)$-modules $D^{\chi_{\delta_0,\mu}}_{k}$ and $W^{\chi_{\delta_0,\mu}}_{k}$ can be described as follows:
		\begin{enumerate}
			\item The module $D^{\chi_{\delta_0,\mu}}_{k}$ has a restriction to the maximal compact subgroup $K = O(2)$
			\[
				D^{\chi_{\delta_0,\mu}}_{k}|_{O(2)} = \bigoplus_{\substack{l\equiv k\text{ }\mathrm{mod}\text{ }2\\l\geq k}} \sigma_{l}.
			\]
			\item The module $W^{\chi_{\delta_0,\mu}}_{k}$ has restriction to the maximal compact subgroup
			\[
				W^{\chi_{\delta_0,\mu}}_{k}|_{O(2)} = \begin{cases}\bigoplus_{\substack{l\equiv k\text{ }\mathrm{mod}\text{ }2\\ 1\leq l\leq k-2}} \sigma_{l} & k\equiv 1\text{ }\mathrm{mod}\text{ }2\\
				\chi_{\delta_0,\mu}|_{O(2)}\oplus\bigoplus_{\substack{l\equiv k\text{ }\mathrm{mod}\text{ }2\\ 2\leq l\leq k-2}} \sigma_{l} & k\equiv 0\text{ }\mathrm{mod}\text{ }2\end{cases}.
			\]
		\end{enumerate}
\end{theorem}
\subsection{An Introduction to Results in this Paper}
In this paper, we mainly deal with the group $SU(2,1)$. The group $SU(2,1)$ is of real rank one, and the decomposition of principal series for $SU(2,1)$ has been studied in \cite{BaldoniEmbeddings}, \cite{WallachComposition2} and \cite{WallachComposition1} . In this paper, I discuss the results for $SU(2,1)$ by passing to the representation theory of their maximal compact subgroups.\\
When the minimal principal series $I(\delta,\lambda)$ of the group $SU(2,1)$ is not irreducible, we have classified 6 families of irreducible sub or quotient $(\frakg,K)$-modules of $I(\delta,\lambda)$ depending on the induction parameters $(\delta,\lambda)$.  In \cite{Kraljevic},\cite{SilvaComposition} and \cite{VoganCharI}, an algorithm to calculate the composition series and the classification of irreducible $(\frakg,K)$-modules for real rank 1 groups like $SU(n,1)$ and $Sp(n,1)$ has already been developed. In this paper, I will utilize the Wigner $D$-functions to perform the calculation explicitly and write down the $K$-types of the irreducible $(\frakg,K)$-modules of $SU(2,1)$ are parametrized by a pair of half integers $(j,n)$ satisfying appropriate parity conditions. These irreducible subquotients are 
\begin{enumerate}
	\item Holomorphic/antiholomorphic discrete series: $V_{\mathrm{disc}\pm}(\delta,\lambda)$,
	\item Quaternionic discrete series $V_{\mathbb{H}}(\delta,\lambda)$,
	\item Finite dimensional representations $V_{\mathrm{fin}}(\delta,\lambda)$,
	\item Two other irreducible $(\frakg,K)$-modules $Q_{\pm}(\delta,\lambda)$.
\end{enumerate}
 Moreover, we can also compute the intertwining operators for the minimal principal series. This paper has developed a computational technique based on combinatorics to calculate the intertwining integrals for $SU(2,1)$. They will be discussed in Section \label{intertwinesection}. For the long intertwining operator of the minimal principal series of $SU(2,1)$, we have another proof of the well-known result from \cite{WallachComposition1}:
\begin{theorem}\label{LongSU21}
	The long intertwining operator $A(w_0,\delta,\lambda)$ acts on each $W^{(j,n)}_{m_1,m_2}$ as a scalar $\left[A(w,\delta,\lambda)\right]_{m_1}$, with a closed form formula:
	\begin{align}
&\left[A(w,\delta,\lambda)\right]_{m_1} = \nonumber\\
	&
	\frac{ \pi^2 2^{-\lambda -1} \Gamma (\lambda)}{\Gamma \left(1-\frac{\lambda -\delta }{2}\right)\Gamma \left(1-\frac{\lambda+\delta}{2}\right)}\frac{ \Gamma \left(j+m_1-\frac{\lambda+\delta }{2}+1\right)\Gamma \left(j-m_1-\frac{\lambda-\delta }{2}+1\right)}{
   \Gamma \left(j+m_1+\frac{\lambda-\delta}{2}+1\right) \Gamma \left(j-m_1+\frac{\lambda+\delta }{2}+1\right)}
\end{align}
\end{theorem}

I would like to thank Stephen D. Miller for leading me to this project and generously sharing his ideas and guidance throughout this project. I would Siddhartha Sahi, James Lepowsky and Doron Zeilberger for useful conversations and instructions. I also would like to thank Oleg Marichev from Wolfram Research of sharing his insights on special functions with me.

\section{Representation of Compact Groups}\label{cpt}
Let $K\subset G$ be a maximal compact subgroup of $G$. We assume $\text{rank } K=\text{rank } G=r$ and there is a compact Cartan subalgebra $\frakt$ for both $\frakk$ and $\frakg$. We denote the complexification of the Lie algebra $\frakk$ as $\frakk_\CC$ and the root system as $\Delta(\frakk_\CC,\frakt_\CC)$, with a choice of simple roots $\mathcal{S} = \{\alpha_1,\ldots,\alpha_r\}$. We introduce the symbol $\hat{K}$  as the set of all irreducible representations of $K$ up to equivalence. By indexing an irreducible representation $(\tau, V_{\tau})$ by its highest weight $\mu$, we can realize the set $\hat{K}$ as a subset of the weight lattice
\[
	X(\frakk_\CC,\frakt_\CC) = \{\mu\in \frakt_\CC^*|\langle\mu,\check{\alpha}_i\rangle\in\ZZ\text{ for all }\alpha_i\in \mathcal{S}\}.
\]
In this section, we will discuss the representation theory of $U(2)$ and $SU(2)$ in detail and we let $K = U(2)$ and $\frakk = \fraku(2)$.
\subsection{Structure of $SU(2)$ and $U(2)$}
The \emph{Pauli matrices} are the generators of the Lie algebra $\fraksl(2,\CC)$:
\begin{align*}
	\sigma_0 = \begin{spmatrix}1&0\\0&1\end{spmatrix},	\sigma_1 =\begin{spmatrix}0&1\\1&0\end{spmatrix},
	\sigma_2 = \begin{spmatrix}0&-\ii\\\ii&0\end{spmatrix},\sigma_3 = \begin{spmatrix}1&0\\0&-1\end{spmatrix}
\end{align*}
They satisfy the commutation relations
\begin{align}
	[\sigma_0,\sigma_i]&=0\\
	[\sigma_i,\sigma_j] &=\sum_{k=1}^3 2\ii\epsilon_{ijk}\sigma_k\text{ if }i,j\neq 0,\label{PauliCommutation}
\end{align}
where the \emph{Levi-Civita symbol} $\epsilon_{ijk}$ takes value $1$ if $(ijk)$ is an even permutation of $(123)$, $-1$ if $(ijk)$ is an odd permutation of $(123)$, and 0 if two or more elements in $\{i,j,k\}$ are equal. We multiply each Pauli matrix by $\ii/2$, and let $\gamma_i=\frac{\ii}{2}\sigma_i$, then $\gamma_i$ are generators of the real Lie algebra $\fraku(2)\subset \fraksl(2,\CC)$:
\begin{align*}
	\gamma_0 = \begin{spmatrix}\frac{\ii}{2}&0\\0&\frac{\ii}{2}\end{spmatrix},	\gamma_1 =\begin{spmatrix}0&\frac{\ii}{2}\\\frac{\ii}{2}&0\end{spmatrix},
	\gamma_2 = \begin{spmatrix}0&\frac{1}{2}\\-\frac{1}{2}&0\end{spmatrix},\gamma_3 = \begin{spmatrix}\frac{\ii}{2}&0\\0&-\frac{\ii}{2}\end{spmatrix}.
\end{align*}
\indent The Pauli matrices $\gamma_0$ and $\gamma_3$ generate a Cartan subalgebra $\mathfrak{t}$ of $\mathfrak{u}(2)$, and the complexified Lie algebra $\frakk_\CC=\fraku(2)\otimes_\RR\CC = \mathfrak{gl}(2,\CC)$ has the following decomposition
\[
	\frakk_\CC = \frakk^-_\CC\oplus\frakk^+_\CC\oplus\frakt_\CC,
\]
where the positive and negative root spaces are generated by $\gamma_1\mp\ii\gamma_2$
\[
	\frakk^-_\CC=\CC(\gamma_1+\ii\gamma_2), \frakk^+_\CC=\CC(\gamma_1-\ii\gamma_2).
\]
The commutator subalgebra $[\frakk_\CC,\frakk_\CC]$ is semisimple and is isomorphic to $\fraksl(2,\CC)$, and $\frakk_\CC = \CC\gamma_0\oplus [\frakk_\CC,\frakk_\CC]$. The Casimir element $\Omega_K$ in the universal enveloping algebra of $U( [\frakk_\CC,\frakk_\CC])$ of the commutator subalgebra  $ [\frakk_\CC,\frakk_\CC]$ takes the form:
\[
	\Omega_K = -2(\gamma_1^2+\gamma_2^2+\gamma_3^2).
\]
The compact group $SU(2)$ is defined by
\[
	SU(2) = \{\begin{spmatrix}\alpha&-\bar{\beta}\\\beta&\bar{\alpha}\end{spmatrix}| \alpha, \beta\in\CC\text{ and }|\alpha|^2+|\beta|^2=1\}\cong \mathbb{S}^3.
\]
We consider the quaternions
\[
	\mathbb{H} = \{q_0+q_1\ii+q_2\mathsf{j}+q_3\mathsf{k}|q_i\in\RR, \ii^2=\mathsf{j}^2=\mathsf{k}^2=\ii\mathsf{j}\mathsf{k} = -1\},
\]
then there exists a group isomorphism between $SU(2)$ and the unit quaternions $\mathbb{H}^*$:
\[
	\begin{matrix}
	SU(2) &\longrightarrow &\mathbb{H}^* = \{q_0+q_1\ii+q_2\mathsf{j}+q_3\mathsf{k}|q_i\in\RR\text{ and }\sum q_i^2 = 1\}\\
	\begin{spmatrix}\alpha&-\bar{\beta}\\\beta&\bar{\alpha}\end{spmatrix} & \mapsto & \mathrm{Re}\alpha +( \mathrm{Im}\alpha) \ii -(\mathrm{Re}\beta) \mathsf{j} +( \mathrm{Im}\beta)\mathsf{k}
	\end{matrix}.
\]
To obtain a rotational coordinate for $SU(2)$, we need to introduce the two variable $\arctan$ function with range $(-\pi,\pi]$:
\[
	\arctan(x,y) = \mathrm{Arg}(x+\ii y)
\]
where $\mathrm{Arg}$ is the principal value of the argument function, taking value in the range $(-\pi,\pi]$. The multiplication by $\mathbb{S}^1\cong\{e^{-\phi \gamma_3}| -\pi <\psi\leq \pi\}$ on the right defines a \emph{Hopf fibration} of the group $SU(2)$:
\[
\begin{matrix}
	\mathbb{S}^1&\longrightarrow&SU(2)&\longrightarrow&\mathbb{CP}^1&\cong&\mathbb{S}^2\\
	&&\begin{spmatrix}\alpha&-\bar{\beta}\\\beta&\bar{\alpha}\end{spmatrix}&&[\alpha:\beta]&&z=\beta/\alpha
\end{matrix},
\]
where we take
\begin{equation}\phi = \arctan(-\mathrm{Im}(\alpha ) \mathrm{Im}(\beta )+\mathrm{Re}(\alpha ) \mathrm{Re}(\beta ),-\mathrm{Im}(\beta ) \mathrm{Re}(\alpha )-\mathrm{Im}(\alpha ) \mathrm{Re}(\beta )).\end{equation}
We can use the $zyz$ Euler angles $(\psi,\theta,\phi)$ to parametrize a generic element of $SU(2)$. The ranges of these angles are
\[
	\begin{matrix}\phi \in(-\pi,\pi],&\theta\in[0,\pi],&\psi\in(-\pi,3\pi].\end{matrix}
\] If we choose the branch for the $\arccos$ function such that its value lies in the range $[0,\pi]$, and let the other two angles
\begin{align}
	\theta &=  \arccos(1-2|\beta|^2),\\
	\psi &=\arctan\left(\beta  \bar{\alpha}+\alpha  \bar{\beta},\mathrm{Re}(2 \alpha  \mathrm{Im}(\beta )-2 \beta  \mathrm{Im}(\alpha ))\right)+\pi(1-\epsilon(\alpha,\beta)),
\end{align}
where
\[
	\epsilon(\alpha,\beta)=\exp \left(-\frac{1}{2} \ii \left(\mathrm{Arg} (\bar{\alpha}  \beta )-2 \mathrm{Arg} (\bar{\alpha} )+\mathrm{Arg} \left(\overline{\alpha\beta}\right)\right)\right),
\]
then $z=\beta/\alpha = e^{\ii \psi}\tan\frac{\theta}{2}$, and the matrix
\begin{equation}
	\mathcal{U}(\psi,\theta,\phi) = e^{-\psi\gamma_3}e^{-\theta\gamma_2}e^{-\phi\gamma_3} = \begin{spmatrix}
	e^{-\frac{\ii}{2}(\phi+\psi)}\cos\frac{\theta}{2}&-e^{\frac{\ii}{2}(\phi-\psi)}\sin\frac{\theta}{2}\\
	e^{\frac{\ii}{2}(-\phi+\psi)}\sin\frac{\theta}{2}&e^{\frac{\ii}{2}(\phi+\psi)}\cos\frac{\theta}{2}
	\end{spmatrix},\label{GenericSU2}
\end{equation} with $\psi,\theta,\phi$ given by the formulas above, parametrizes a generic element $\begin{spmatrix}\alpha&-\bar{\beta}\\\beta&\bar{\alpha}\end{spmatrix}$ of the group $SU(2)$ with the Euler angles $(\psi,\theta,\phi)$. If we write $\begin{spmatrix}\alpha&-\bar{\beta}\\\beta&\bar{\alpha}\end{spmatrix}$ as a unit quaternion, then the matrix with entries in $\alpha$ and $\beta$ corresponds to the unit quaternion $q_0+q_1\ii+q_2\mathsf{j}+q_3\mathsf{k}$, where
\[
	\begin{matrix}q_0=\cos \left(\frac{\theta }{2}\right) \cos \left(\frac{\psi +\phi }{2}\right),&q_1=-\cos \left(\frac{\theta }{2}\right) \sin
   \left(\frac{\psi +\phi }{2}\right),\\
   q_2=-\sin \left(\frac{\theta }{2}\right) \cos \left(\frac{\phi -\psi }{2}\right),&q_3=\sin
   \left(\frac{\theta }{2}\right) \sin \left(\frac{-\phi +\psi }{2}\right).
   \end{matrix}
\]The center of the Lie group $U(2)$ is isomorphic to $U(1)$, which can be parametrized by the exponentiation $e^{-\zeta\gamma_0}$ of the matrix $\gamma_0$. If we multiply $e^{-\zeta\gamma_0}$ to the matrix $\mathcal{U}(\psi,\theta,\phi)$,  we obtain a generic element of the group $U(2)$:
\begin{align}
	e^{-\zeta\gamma_0}\mathcal{U}(\psi, \theta, \phi)=\begin{spmatrix}
	e^{\frac{\ii}{2}(-\zeta-\phi-\psi)}\cos\frac{\theta}{2}&-e^{\frac{\ii}{2}(-\zeta+\phi-\psi)}\sin\frac{\theta}{2}\\
	e^{\frac{\ii}{2}(-\zeta-\phi+\psi)}\sin\frac{\theta}{2}&e^{\frac{\ii}{2}(-\zeta+\phi+\psi)}\cos\frac{\theta}{2}
	\end{spmatrix}\label{GenericSU2Matrix}.
\end{align}
\subsection{Realization of Irreducible Representations}\label{Irrep}
The standard representation $(\pi,W)$ of $U(2)$ can be realized in the space of linear polynomials in two variables. To be precise, let $z$ be the column vector of two variables $(z_1,z_2)^t$, with $\pi$ acting on any linear polynomial $f(z)$ via:
\[
	\pi(g)f(z)=f(g^{-1}z).
\]
For $k\in\ZZ$, denote by $\det^k$ the one dimensional representation on which $g\in U(2)$ acts by scalar multiplication with $(\det g)^{k}$. Letting $j\in\frac{1}{2}\mathbb{N}$ and $n\in\frac{1}{2}\ZZ$ with $j+n\in \ZZ$, an arbitrary irreducible representation $\pi_{j,n}$ of $U(2)$ can be realized on the space $\Sym^{2j}W\otimes \det^{(j+n)}$, which is isomorphic to $\Sym^{2j} W$ as a vector space. $U(2)$ acts by right regular action on any degree $2j$ homogeneous polynomial $f\in \Sym^{2j}(W)\otimes \det^{(j+n)}$ in 2 variables $z_1,z_2$:
\[
	\pi_{j,n}(g)f(z)=(\det g)^{j+n}f(g^{-1}z).
\]
Let $m\in\frac{1}{2}\ZZ$ such that $-j\leq m\leq j$ and $j\pm m$ are integers, the weight basis $\{v^j_m\}_{-j\leq m\leq j}$ for $\Sym^{2j}W \otimes \det^{(j+n)}$ is defined as:
\[
	v^j_m = \frac{z_1^{j-m}z_2^{j+m}}{\sqrt{(j-m)!(j+m)!}}.
\]
The Lie algebra $\mathfrak{u}(2)$ acts on the weight basis $\{v^j_m\}_{-j\leq m\leq j}$ of $\Sym^{2j}W\otimes \det^{(j+n)}$ as linear operators:
\begin{align}
	\gamma_0 v^{j}_m &= \ii n v^{j}_m\label{normcond1}\\
	(\gamma_1\pm \ii \gamma_2) v^{j}_m &= -\ii\sqrt{(j\mp m)(j\pm m+1)} v^{j}_{m\pm1}\label{normcond2}\\
	\gamma_3  v^{j}_m&= \ii m  v^{j}_m\label{normcond3}.
\end{align} We define a hermitian inner product $\langle v^j_{m_1}, v^j_{m_2} \rangle = \delta_{m_1,m_2}
$ on $\Sym^{2j} W\otimes \det^{(j+n)}$, and we require the inner product to be linear in the first argument, and conjugate linear in the second argument. We notice that $\gamma_0$ and $\gamma_3$ acts on the weight vectors $v^j_{m}$ of $\Sym^{2j}W \otimes\det^{(j+n)}$ by multiplication of a purely imaginary number $\ii m$. Moreover, from (\ref{normcond2}), we can see that the matrices of $\gamma_1$ and $\gamma_2$ actions are unitary under this hermitian inner product:
\begin{align}
	\langle \gamma_1 v^j_{m_1},&v^j_{m_2}\rangle =\nonumber\\
	& \frac{1}{2\ii}\left(\sqrt{(j-m_1)(j+m_1+1)}\delta_{m_1+1,m_2}+\sqrt{(j-m_2)(j+m_2+1)}\delta_{m_1,m_2+1}\right)\\
	\langle \gamma_2 v^j_{m_1},&v^j_{m_2}\rangle =\nonumber\\
	& \frac{1}{2}\left(-\sqrt{(j-m_1)(j+m_1+1)}\delta_{m_1+1,m_2}+\sqrt{(j-m_2)(j+m_2+1)}\delta_{m_1,m_2+1}\right).
\end{align}
Therefore, the action of the Lie algebra elements $\gamma_i\in\fraku(2)$ on $\Sym^{2j} W\otimes  \det^{(j+n)}$ is unitary. Using this inner product, we can define the \emph{Wigner $D$-functions}  $W^{(j,n)}_{m_1,m_2}$ as the matrix coefficients of the irreducible representation $\pi_{j,n}$:
\begin{align}
	W^{(j,n)}_{m_1,m_2}(\zeta,\psi,\theta,\phi)&=\langle v^j_{m_1},e^{-\gamma_0\zeta}\mathcal{U}(\psi,\theta,\phi)v^j_{m_2}\rangle\nonumber\\
	&= c^j_{m_1}c^j_{m_2}e^{\ii n\zeta}e^{\ii(m_1\psi+m_2\phi)}d^{(j,n)}_{m_1,m_2}(\theta)\label{WignerDefinition},
\end{align}
where $c^j_m = \sqrt{(j+m)!(j-m)!}$ is a normalization factor, and the function $d^{(j,n)}_{m_1,m_2}(\theta)$ is given by the trigonometric polynomial
\begin{align}
d^{(j,n)}_{m_1,m_2}(\theta)&=\sum_{p=\max(0,m_1-m_2)}^{\min(j-m_2,j+m_1)}\frac{(-1)^{m_2-m_1+p}}{(j+m_1-p)!p!(m_2-m_1+p)!(j-m_2-p)!}\nonumber\\&\sin^{m_2-m_1+2p}\left(\frac{\theta}{2}\right)\cos^{2 j + m_1 - m_2 - 2 p}\left(\frac{\theta}{2}\right).\label{DPart}
\end{align}
The Wigner $D$-function $W^{(j,n)}_{m_1,m_2}$ satisfies the following properties:
\begin{enumerate}
\item {\bf Jacobi polynomials and $d^{(j,n)}_{m_1,m_2}(\theta)$ }\\
The sum $d^{(j,n)}_{m_1,m_2}(\theta)$ has an expression in terms of the hypergeometric function $_2F_1$:
\begin{align}
	&d^{(j,n)}_{m_1,m_2}(\theta) =\nonumber\\
	&\begin{cases}\frac{\sin^{m_1-m_2}\left(\frac{\theta}{2}\right)\cos^{2j-m_1+m_2}\left(\frac{\theta}{2}\right)}{(j-m_1)!(m_1-m_2)!(j+m_2)!}{}_2F_1\left(\begin{smallmatrix}-j+m_1,-j-m_2\\1+m_1-m_2\end{smallmatrix};-\tan^2\left(\frac{\theta}{2}\right)\right)&m_1>m_2\\
	\frac{(-1)^{-m_1+m_2}\sin^{-m_1+m_2}\left(\frac{\theta}{2}\right)\cos^{2j+m_1-m_2}\left(\frac{\theta}{2}\right)}{(j+m_1)!(-m_1+m_2)!(j-m_2)!}{}_2F_1\left(\begin{smallmatrix}-j-m_1,-j+m_2\\1-m_1+m_2\end{smallmatrix};-\tan^2\left(\frac{\theta}{2}\right)\right)&m_1\leq m_2\end{cases}\label{dashyper}
\end{align}
For $n\geq 0$ and for $\alpha,\beta\in\RR$, the \emph{Jacobi polynomials} $P_{n}^{{(\alpha ,\beta )}}(z)$ are a class of orthogonal polynomials defined in \cite{AbStegun} as
\begin{align}
	P_{n}^{{(\alpha ,\beta )}}(z)={\frac {\Gamma (\alpha +n+1)}{n!\,\Gamma (\alpha +\beta +n+1)}}\sum _{{m=0}}^{n}{n \choose m}{\frac {\Gamma (\alpha +\beta +n+m+1)}{\Gamma (\alpha +m+1)}}\left({\frac {z-1}{2}}\right)^{m}.\label{Jacobi1}
\end{align}
In \cite{WolframAlpha}, the Jacobi polynomial is also defined as
\begin{align}\label{Jacobi2}
P^{\alpha,\beta}_n(x)&=\binom{n+\alpha}{n}\left(\frac{x+1}{2}\right)^n {}_2F_1\left(\begin{smallmatrix}-n,-n-\beta,\alpha+1\end{smallmatrix};\frac{x-1}{x+1}\right).
\end{align}
These two definitions of Jacobi polynomials are equivalent. By the definition of Jacobi polynomials and the transformation rules of hypergeometric functions, the function $d^{(j,n)}_{m_1,m_2}(\theta)$ can also be expressed in terms of Jacobi polynomials:
\begin{align}
	&d^{(j,n)}_{m_1,m_2}(\theta)=\frac{\left(\sin\frac{\theta}{2}\right)^{m_1-m_2}\left(\cos\frac{\theta}{2}\right)^{m_1+m_2}}{(j+m_2)!(j-m_2)!}P^{(m_1-m_2,m_1+m_2)}_{j-m_1}(\cos\theta).\label{WignerJacobi}
\end{align}
%They are orthogonal with respect to the measure $(1-x)^\alpha(1+x)^\beta\dd x$ on the interval $[-1,1]$:\cite{AbStegun}
%\begin{align*} \int _{-1}^{1}(1-x)^{\alpha }(1+x)^{\beta }&P_{m}^{(\alpha ,\beta )}(x)P_{n}^{(\alpha ,\beta )}(x)\,dx=\\
%&{\frac {2^{\alpha +\beta +1}}{2n+\alpha +\beta +1}}{\frac {\Gamma (n+\alpha +1)\Gamma (n+\beta +1)}{\Gamma (n+\alpha +\beta +1)n!}}\delta _{nm}.\end{align*}

\item {\bf Multiplicativity}
If an element $k$ of $U(2)$ is expressed in terms of Euler angles as $k = e^{-\gamma_0\zeta}\mathcal{U}(\psi,\theta,\phi)$, we can replace the notation $W^{(j,n)}_{m_1,m_2}(\zeta,\psi,\theta,\phi)$ by $W^{(j,n)}_{m_1,m_2}(k)$. For any $k_1,k_2\in K$, since $W^{(j,n)}_{m_1,m_2}$ are the matrix coefficients of the representation $\pi_{j,n}$, we can use the multiplicative property of matrix coefficients to write $W^{(j,n)}_{m_1,m_2}(k_1k_2)$ as a sum:
\begin{align*}
	W^{(j,n)}_{m_1,m_2}(k_1k_2)&=\langle v^j_{m_1}, k_1 k_2v^j_{m_2}\rangle\\
	&=\sum_{\substack{-j\leq m_3\leq j\\j+m_3\in\ZZ}}\langle v^j_{m_1},k_1 v^j_{m_3}\rangle\langle v^j_{m_3}, k_2v^j_{m_2}\rangle\\
	&=\sum_{\substack{-j\leq m_3\leq j\\j+m_3\in\ZZ}} W^{(j,n)}_{m_1,m_3}(k_1) W^{(j,n)}_{m_3,m_2}(k_2).
\end{align*}

\item {\bf Inverse Matrix}\\
Recall that the basis $\{v^j_{m}\}$ and the inner product on $\pi_{j,n}$ are chosen so that the action of $U(2)$ is unitary. Therefore, the Wigner $D$-functions satisfies the unitarity property:
\begin{equation}
\sum_{\substack{-j\leq m_3\leq j\\j+m_3\in\ZZ}}W^{(j,n)}_{m_1,m_3}(k)\overline{W^{(j,n)}_{m_2,m_3}}(k) = \delta_{m_1,m_2}.\label{Unitarity1}
\end{equation}
This relation is equivalent to:
\[
	\overline{W^{(j,n)}_{m_2,m_1}(k)}=W^{(j,n)}_{m_1,m_2}(k^{-1}).
\]
To see such transformation rule directly, we can switch $m_1$ and $m_2$ in the expression (\ref{dashyper}) of $d^{(j,n)}_{m_1,m_2}(\theta)$ in terms of hypergeometric functions, and obtain the formula\[d^{(j,n)}_{m_1,m_2}(\theta)=(-1)^{m_2-m_1}d^{(j,n)}_{m_2,m_1}(\theta).\]
It follows from the previous formula and the definition of Wigner $D$-functions that
\begin{equation}
	 (-1)^{m_2-m_1}W^{(j,-n)}_{-m_1,-m_2}(k)=W^{(j,n)}_{m_2,m_1}(k^{-1}).\label{UnitarityWigner}
\end{equation}
which is equivalent to the formula (\ref{Unitarity1}) above.

\item {\bf Differential Equations}\\
The right and left regular action $r(k)$ and $l(k)$ by $K$ on any function $f\in C^\infty(K)$ are given by
\begin{align*}
	\text{Right action: }(r(k)f)(g) &= f(gk)\\
	\text{Left action: }(l(k)f)(g) &= f(k^{-1}g).
\end{align*}
The corresponding action of the Lie algebra $\fraku(2)$ as differential operators are denoted by $\dd l$ and $\dd r$, and we can extend the action to $\fraksl(2,\CC)=\fraku(2)\otimes\CC$ by making $\dd l$ linear and $\dd r$ linear under multiplication by scalars:
\[
	\begin{matrix}\dd l(\alpha X) = \alpha \dd l(X),& \dd r(\alpha X) = \alpha \dd r(X).\end{matrix}
\] The differential operators can be written down explicitly using the Euler angle coordinates. Then by direct calculation, these differential operators are
\begin{align*}
	&\begin{matrix}\dr(\gamma_0)=-\frac{\partial}{\partial\zeta}& 	\dl(\gamma_0)=\frac{\partial}{\partial\zeta}\\
	\dr(\gamma_3)=-\frac{\partial}{\partial\phi}&	\dl(\gamma_3)=\frac{\partial}{\partial\psi}\end{matrix}\\
	&\dr(\gamma_1\pm \ii \gamma_2)=e^{\mp \ii\phi}(-\cot\theta\frac{\partial}{\partial\phi}\mp\ii \frac{\partial}{\partial\theta}+\csc\theta\frac{\partial}{\partial\psi})\\
	&\dl(\gamma_1\pm \ii \gamma_2)=e^{\pm \ii\psi}(\csc\theta \frac{\partial}{\partial\phi}\pm \ii \frac{\partial}{\partial\theta}-\cot\theta\frac{\partial}{\partial\psi}).
\end{align*}
Comparing with the action of $\gamma_i$ on the weight basis $v^j_m$ in (\ref{normcond1})-(\ref{normcond3}), the action of the differential operators $\dr(\gamma_i)$ and $\dl(\gamma_i)$ on the Wigner $D$-functions are
\begin{align}
\begin{matrix}
\dr(\gamma_0)W^{(j,n)}_{m_1,m_2}=-\ii nW^{(j,n)}_{m_1,m_2}& \dl(\gamma_0)W^{(j,n)}_{m_1,m_2}=\ii nW^{(j,n)}_{m_1,m_2}\\
	\dr(\gamma_3)W^{(j,n)}_{m_1,m_2}=-\ii m_2W^{(j,n)}_{m_1,m_2}& \dl(\gamma_3)W^{(j,n)}_{m_1,m_2}=\ii m_1W^{(j,n)}_{m_1,m_2}\\
\end{matrix}\label{CompactAction:1}\\
	\dr(\gamma_1\pm \ii\gamma_2)W^{(j,n)}_{m_1,m_2}=\ii\sqrt{(j\pm m_2)(j\mp m_2+1)}W^{(j,n)}_{m_1,m_2\mp 1}\label{CompactAction:2}\\
	\dl(\gamma_1\pm \ii\gamma_2)W^{(j,n)}_{m_1,m_2}=-\ii\sqrt{(j\mp m_1)(j\pm m_1+1)}W^{(j,n)}_{m_1\pm 1,m_2}\label{CompactAction:3}.
\end{align}
\item {\bf Basis for $L^2(K)$}\\
By the Peter-Weyl Theorem for $K=U(2)$, as the matrix coefficients for the finite dimensional representations of $U(2)$, the Wigner $D$-functions $W^{(j,n)}_{m_1,m_2}$ provides a Hilbert space basis for $L^2(K)$:
\begin{align}
	L^2(K) =\widehat{\bigoplus}_{\substack{j\in\frac{1}{2}\ZZ_{\geq 0}\\n\in\frac{1}{2}\ZZ, j+n\in\ZZ\\m_1,m_2\in\{-j,-j+1,\ldots,j\}}}\CC W^{(j,n)}_{m_1,m_2}.\label{HilbertBasis}
\end{align}

\end{enumerate}

\subsection{Tensor products and Clebsch-Gordan coefficients}\label{CGCoef}
For any two irreducible representations $V^{j_1}\cong \Sym^{2j_1}W$ and $V^{j_2}\cong\Sym^{2j_2} W$ of $SU(2)$, we choose the weight basis $\{v^{j_1}_{m_1}\}$ and $\{v^{j_2}_{m_2}\}$ of the two spaces, properly normalized such that $\gamma_i$ acts as in (\ref{normcond1})-(\ref{normcond3}). The tensor product $V^{j_1}\otimes V^{j_2}$ has two sets of basis: the \emph{pure tensors} $v^{j_1}_{m_1}\otimes v^{j_2}_{m_2}$, and the weight basis of the irreducible constituents $V^{J}$ of the tensor product $V^{j_1}\otimes V^{j_2}$. In fact, for the compact group $SU(2)$, each irreducible constituent in the decomposition
\[
	V^{j_1}\otimes V^{j_1} = \bigoplus_{J} (V^{J})^{\oplus m_J}
\]
has multiplicity $m_J = 0$ or $1$, with $m_J=1$ if and only if $J$ satisfies the following two properties:
\begin{enumerate}
	\item $|j_1-j_2|\leq J \leq j_1+j_2$
	\item $J-|j_1-j_2|\in\ZZ$.
\end{enumerate}
To be more precise about the relationship between the pure tensor basis and the weight basis, we can expand the pure tensor basis into a linear combination of the weight basis:
\[
	v^{j_1}_{m_1}\otimes v^{j_2}_{m_2} \equiv \sum_{\substack{|j_1-j_2|\leq J \leq j_1+j_2\\J-|j_1-j_2|\in\ZZ}}\left( \begin{smallmatrix}J,M\\j_1,m_1,j_2,m_2\end{smallmatrix}\right) v^J_{M}.
\]
The coefficient $\left( \begin{smallmatrix}J,M\\j_1,m_1,j_2,m_2\end{smallmatrix}\right) $ in the expression above is called the \emph{Clebsch-Gordan coefficient}. It is zero except when $M = m_1+m_2$. Moreover, the Clebsch-Gordan coefficients can also be used to write the product of Wigner $D$-functions as a linear combinations of Wigner $D$-functions. If we introduce the inner product on the tensor product $V^{j_1}\otimes V^{j_2}$ such that
\[
\langle v^{j_1}_{m_{12}}\otimes v^{j_2}_{m_{22}}, v^{j_1}_{m_{11}}\otimes v^{j_2}_{m_{21}} \rangle= \langle v^{j_1}_{m_{12}}, v^{j_1}_{m_{11}}\rangle\langle v^{j_2}_{m_{22}}, v^{j_2}_{m_{21}}\rangle,
\] the product $W^{(j_1,n_1)}_{m_{11},m_{12}}W^{(j_2,n_2)}_{m_{21},m_{22}}$ of Wigner $D$-functions is thus a matrix coefficient of the representation $V^{j_1}\otimes V^{j_2}$:
\begin{align*}
	W^{(j_1,n_1)}_{m_{11},m_{12}}W^{(j_2,n_2)}_{m_{21},m_{22}} &= \langle v^{j_1}_{m_{12}}, k v^{j_1}_{m_{11}}\rangle\langle v^{j_2}_{m_{22}}, k v^{j_2}_{m_{21}}\rangle\\
	&=\langle v^{j_1}_{m_{12}}\otimes v^{j_2}_{m_{22}}, k (v^{j_1}_{m_{11}}\otimes v^{j_2}_{m_{21}})\rangle.
\end{align*}
Combining this with the expansion of pure tensors into linear combinations of weight basis:
\begin{align*}
	v^{j_1}_{m_{11}}\otimes v^{j_2}_{m_{21}} &= \sum_{J_1}\left( \begin{smallmatrix}J_1, M_1\\j_1,m_{11},j_2,m_{21}\end{smallmatrix}\right) v^{J_1}_{M_1}\\
	v^{j_1}_{m_{21}}\otimes v^{j_2}_{m_{22}} &= \sum_{J_2}\left( \begin{smallmatrix}J_2, M_2\\j_1,m_{12},j_2,m_{22}\end{smallmatrix}\right) v^{J_2}_{M_2},
\end{align*}
and since the matrix coefficients of nonisomorphic irreducible representations are orthogonal to each other, the product of Wigner $D$-functions $W^{(j_1,n_1)}_{m_{11},m_{12}}W^{(j_2,n_2)}_{m_{21},m_{22}}$ can be written as
\begin{align}
	W^{(j_1,n_1)}_{m_{11},m_{12}}W^{(j_2,n_2)}_{m_{21},m_{22}} = \sum_{\substack{|j_1-j_2|\leq J \leq j_1+j_2\\J-|j_1-j_2|\in\ZZ\\M_1=m_{11}+m_{21}\\M_2=m_{12}+m_{22}}}\left( \begin{smallmatrix}J, M_1\\j_1,m_{11},j_2,m_{21}\end{smallmatrix}\right) \left( \begin{smallmatrix}J, M_2\\j_1,m_{12},j_2,m_{22}\end{smallmatrix}\right)W^{(J,n_1+n_2)}_{M_1,M_2}.\label{ProductClebschGordan}
\end{align}
The Clebsch-Gordan coefficient $\left( \begin{smallmatrix}J, M\\j_1,m_1,j_2,m_2\end{smallmatrix}\right)$ can be related to the \emph{Wigner $3j$-symbols} $\left(\begin{smallmatrix}j_1&j_2&J\\m_1&m_2&-M\end{smallmatrix}\right)$ in the following way:
\[
	\left(\begin{smallmatrix}j_1&j_2&J\\m_1&m_2&-M\end{smallmatrix}\right) = \frac{(-1)^{j_2-j_1-M}}{\sqrt{2 J+1}} \left( \begin{smallmatrix}J, M\\j_1,m_1,j_2,m_2\end{smallmatrix}\right).
\]
The Wigner $3j$-symbol $\left(\begin{smallmatrix}j_1&j_2&j_3\\m_1&m_2&m_3\end{smallmatrix}\right)\neq 0$ is nonzero if and only if the following conditions are satisfied:
\begin{enumerate}
	\item $m_i = -j_i, -j_i+1, \ldots, j_i-1, j_i$;
	\item $m_1+m_2+m_3=0$;
	\item $|j_1-j_2|\leq j_3\leq j_1+j_2$;
	\item $j_1+j_2+j_3\in\ZZ$.
\end{enumerate}
There is also a recursion relation of the Wigner $3j$-symbols, written in a symmetric manner as
\begin{align*}
	&\sqrt{(j_1\mp m_1)(j_1\pm m_1+1)}\left(\begin{smallmatrix}j_1&j_2&j_3\\m_1\pm 1&m_2&m_3\end{smallmatrix}\right)+\sqrt{(j_2\mp m_2)(j_2\pm m_2+1)}\left(\begin{smallmatrix}j_1&j_2&j_3\\m_1&m_2\pm 1&m_3\end{smallmatrix}\right)\\&+\sqrt{(j_3\mp m_3)(j_3\pm m_3+1)}\left(\begin{smallmatrix}j_1&j_2&j_3\\m_1&m_2&m_3\pm 1\end{smallmatrix}\right)=0.
\end{align*}
%The recursion starts from cases when $m_1, m_2=\pm j$:
%\begin{align*}
%	\sqrt{(j_2\mp m_2)(j_2\pm m_2+1)}\begin{spmatrix} j_1&j_2&j_3\\j_1&m_2\pm 1&-j_1-m_2\mp 1\end{spmatrix}+\sqrt{(j_3\mp m_3)(j_3\pm m_3+1)}\begin{spmatrix} j_1&j_2&j_3\\j_1&m_2&-j_1-m_2\end{spmatrix}
%\end{align*}
The Clebsch-Gordan coefficients $\begin{spmatrix}j+j_0,m_1+m_2\\j,m_1,\frac{1}{2},m_2\end{spmatrix}$ and $\begin{spmatrix}j+j_0,m_1+m_2\\j,m_1,1,m_2\end{spmatrix}$ are listed in Table \ref{TableCG1} and Table \ref{TableCG2}.\\
\begin{table}[h]
	\centering
	\begin{tabular}{ccc}
		$\begin{spmatrix}j+j_0,m_1+m_2\\j,m_1,\frac{1}{2},m_2\end{spmatrix}$& $m_2=-\frac{1}{2}$ & $m_2=+\frac{1}{2}$\\\hline
		$j_0=-\frac{1}{2}$ & $\sqrt{\frac{j+m_1}{2j+1}}$ & $-\sqrt{\frac{j-m_1}{2j+1}}$ \\
		$j_0 = \frac{1}{2}$ & $\sqrt{\frac{j-m_1+1}{2j+1}}$ &  $\sqrt{\frac{j+m_1+1}{2j+1}}$
	\end{tabular}
	\caption{Table for Clebsch-Gordan coefficients of $V^{j}\otimes V^{\frac{1}{2}}$}\label{TableCG1}
\end{table}

\begin{table}[h]
\centering
	\begin{tabular}{cccc}
		$\begin{spmatrix}j+j_0,m_1+m_2\\j,m_1,1,m_2\end{spmatrix}$& $m_2=-1$ & $m_2=0$ & $m_2=1$\\\hline
		$j_0=-1$ & $\sqrt{\frac{(j+m_1)(j+m_1-1)}{2j(2j+1)}}$ & $-\sqrt{\frac{(j-m_1)(j+m_1)}{j(2j+1)}}$ & $\sqrt{\frac{(j-m_1)(j-m_1-1)}{2j(2j+1)}}$\\
		$j_0=0$ & $\sqrt{\frac{(j+m_1)(j-m_1+1)}{2j(j+1)}}$ & $\frac{m_1}{\sqrt{j(j+1)}}$ & $-\sqrt{\frac{(j-m_1)(j+m_1+1)}{2j(j+1)}}$\\
		$j_0 =1$ & $\sqrt{\frac{(j-m_1+1)(j-m_1+2)}{(2j+2)(2j+1)}}$ &  $\sqrt{\frac{(j-m_1+1)(j+m_1+1)}{(j+1)(2j+1)}}$ & $\sqrt{\frac{(j+m_1+1)(j+m_1+2)}{(2j+2)(2j+1)}}$
	\end{tabular}.
	\caption{Table for Clebsch-Gordan coefficients of $V^{j}\otimes V^{1}$}\label{TableCG2}
\end{table}

\section{The Group $SU(2,1)$}
The group $SU(2,1)$ is a real form of $SL(3,\CC)$. It is the fixed point of the antiholomorphic involution $\sigma:g\mapsto J (\bar{g}^t)^{-1} J^{-1}$ in $SL(3,\CC)$:
\[
	SU(2,1)=\{g\in SL(3,\CC)|\bar{g}^t J g=J\}
\]
where $J$ is the matrix
\[
	J=\diag(1,1,-1).
\]
The real Lie algebra $\frakg$ of $G=SU(2,1)$ is
\[
	\frakg=\fraksu(2,1)=\{X\in \fraksl(3,\CC)| \bar{X}^t J+JX=0\}.
\]
The structure theory and representation theory of $SU(2,1)$ has been discussed in \cite{ShahidiEisenstein}.
\subsection{The structure of $SU(2,1)$}
First we consider the complex Lie algebra $\fraksl(3,\CC)$ and choose the following data
\begin{enumerate}
	\item A Cartan subalgebra $\frakh_\CC$ generated by $H_{\alpha_1} = E_{11}-E_{22}$ and $H_{\alpha_2} = E_{22}-E_{33}$;
	\item The fundamental weights $\varpi_1,\varpi_2$ in $\frakh_\CC^*$ as dual basis for $H_{\alpha_1},H_{\alpha_2}$, satisfying $\langle \varpi_i,H_{\alpha_j}\rangle = \delta_{ij}$ with $i,j\in\{1,2\}$;
	\item The simple roots $\alpha_1 = 2\varpi_1-\varpi_2, \alpha_2 = -\varpi_1+2\varpi_2$;
	\item The set of positive roots $\Delta^+(\frakg_\CC,\frakh_\CC) = \{\alpha_1,\alpha_2, \alpha_1+\alpha_2\}$;
	\item $\rho_{\CC} = \frac{1}{2}\sum_{\alpha\in\Delta^+(\frakg_\CC,\frakh_\CC)}\alpha = \alpha_1+\alpha_2$;
	\item The generators for the positive root spaces $X_{\alpha_1} =\begin{spmatrix}0&1&0\\0&0&0\\0&0&0\end{spmatrix}$, $X_{\alpha_2} =\begin{spmatrix}0&0&0\\0&0&1\\0&0&0\end{spmatrix}$, $X_{\alpha_1+\alpha_2}=[X_{\alpha_1},X_{\alpha_2}]=\begin{spmatrix}0&0&1\\0&0&0\\0&0&0\end{spmatrix}$
	\item Choose the generators of the negative root spaces $X_{-\alpha} = X_{\alpha}^t$.
\end{enumerate}

The maximal compact subgroup $K$ is defined as the set of fixed points of the Cartan involution $\theta: g\mapsto -\bar{g}^t$ on $G$:
\begin{align*}
	K&=G^\theta=(U(2)\times U(1))/U(1)=\{\begin{spmatrix}A&0\\0&\det(A)^{-1}\end{spmatrix}|A\in U(2)\}\cong U(2).
\end{align*}
The induced isomorphism from the Lie algebra $\fraku(2)$ to $\frakk$ sends the Pauli matrices $\gamma_i$ to the following elements in $\frakk$:
\begin{align*}
	U_0=\frac{1}{2}\begin{spmatrix}\ii&0&0\\0&\ii&0\\0&0&-2\ii\end{spmatrix},
	U_1=\frac{1}{2}\begin{spmatrix}0&\ii&0\\\ii&0&0\\0&0&0\end{spmatrix},
	U_2=\frac{1}{2}\begin{spmatrix}0&1&0\\-1&0&0\\0&0&0\end{spmatrix},
	U_3=\frac{1}{2}\begin{spmatrix}\ii&0&0\\0&-\ii&0\\0&0&0\end{spmatrix}.
\end{align*}
Under this isomorphism, the generators of the compact Cartan subalgebra $\frakt\subset \fraksu(2,1)$ are $\ii H_{\alpha_1} = 2 U_3$ and $\ii H_{\alpha_2} = U_0-U_3$. The complexification of the compact Cartan subalgebra $\frakt_\CC$ is thus exactly the same as the Cartan subalgebra $\frakh_\CC$ of $\fraksl(3,\CC)$ defined above.\\

The -1 eigenspace $\frakp$ of $\theta$ in the Cartan decomposition $\frakg = \frakk\oplus\frakp$ can be described explicitly as the following space of $3\times 3$ matrices with complex entries:
\[
	\frakp = \left\{\begin{spmatrix}0&0&z_1\\0&0&z_2\\\bar{z}_1&\bar{z}_2&0\end{spmatrix}|z_1,z_2\in \CC\right\}.
\]
The matrix $\begin{spmatrix}0&0&z_1\\0&0&z_2\\\bar{z}_1&\bar{z}_2&0\end{spmatrix}$ transforms under the adjoint action by an element $e^{a U_0+b U_3}$ in the Cartan subgroup $T\subset K$,  
\begin{equation}
	\Ad(e^{a U_0+b U_3})\begin{spmatrix}0&0&z_1\\0&0&z_2\\\bar{z}_1&\bar{z}_2&0\end{spmatrix}=\begin{spmatrix}0&0&e^{\ii\frac{3a+b}{2}}z_1\\0&0&e^{\ii\frac{3a-b}{2}}z_2\\e^{-\ii\frac{3a+b}{2}}\bar{z}_1&e^{-\ii\frac{3a-b}{2}}\bar{z}_2&0\end{spmatrix}.\label{Z1Z2}
\end{equation}
Therefore, under the adjoint action of $K\cong U(2)$, we can consider the complexified space $\frakp_\CC = \frakp\otimes\CC$ as a 4-dimensional representation of $K$. Under the usual pairing between $\frakt_\CC$ and $\frakt_\CC^*$, the value of any character $\chi = \chi_1\varpi_1+\chi_2\varpi_2\in \frakt_\CC^*$ on $\frakt$ is given by \[\chi(a U_0+b U_3) = \ii \frac{a+b}{2} \chi_1+ \ii a \chi_2.\]  From this formula, the two simple roots $\alpha_1,\alpha_2\in \Delta(\frakg_\CC,\frakt_\CC)$ act on $aU_0+b U_3$ by
\begin{align}
	\alpha_1(a U_0+b U_3) &= \ii b,\label{Beta1}\\
	\alpha_2(a U_0+b U_3) &= \ii\frac{3a-b}{2}\label{Beta2}
\end{align}
and the action of the highest root is
\begin{align}
	(\alpha_1+\alpha_2)(a U_0+b U_3)&=\ii\frac{3a+b}{2}.\label{Beta1Beta2}
\end{align}

\subsection{The Cartan Subgroups of $SU(2,1)$}
There are two conjugacy classes of Cartan subgroups of $SU(2,1)$: the compact Cartan subgroup isomorphic to $U(1)\times U(1)$, and the maximally noncompact Cartan subgroup isomorphic to $U(1)\times \RR^{\times}$.

\subsubsection{The maximally compact Cartan subgroup}
In the root system $\Delta(\frakg_\CC,\frakt_\CC)$ with respect to the Cartan subalgebra $\frakt_\CC$, the Vogan diagram of $SU(2,1)$
\begin{center}
\begin{tikzpicture}
	\draw
	(0,0) circle [radius=.1] node [above] {$\alpha_1$};
	\draw[fill=black]
	(1,0) circle [radius=.1] node [above] {$\alpha_2$};
	\draw
	(0.1,0) --++ (1,0);
%	(0.1,+.04) --++ (1,0);
%	\draw
%	(.5,0) --++ (60:.2)
%	(.5,0) --++ (-60:.2);
\end{tikzpicture}
\end{center}
specifies an imaginary compact simple root $\alpha_1$ and an imaginary noncompact simple root $\alpha_2$. Since the root vectors of $\alpha_1+\alpha_2$ can always be written as the commutators of root vectors of $\alpha_1$ and $\alpha_2$, the root $\alpha_1+\alpha_2$ is a noncompact root. We have thus obtained the set of positive compact roots $\Delta^+_{c}$ and noncompact roots $\Delta^+_{nc}$:
\begin{align*}
	\Delta^+_{c}(\frakg_\CC,\frakh_\CC) &= \{\alpha_1\}\\
	\Delta^+_{nc} (\frakg_\CC,\frakh_\CC)&= \{\alpha_2,\alpha_1+\alpha_2\}.
\end{align*}
Comparing with the discussion on $U(2)$-representations in Section {\ref{cpt}}, where a weight $m\in\frac{1}{2}\ZZ$ of $\frakp_\CC$ is equal to the half integer $-\ii\alpha(U_3)$, and the central $U(1)$-character $n$ is equal to $-\ii\alpha(U_0)$. Based (\ref{Beta1})-(\ref{Beta1Beta2}), any noncompact imaginary root \[\alpha \in\{ \pm\alpha_2, \pm\alpha_1\pm\alpha_2\}\] correspond to a pair of half integers denoted by $(m_\alpha,n_\alpha) = (-\ii\alpha(U_3),-\ii\alpha(U_0))$. Therefore, we have established the correspondence between noncompact imaginary roots with pairs of half integers $(m_\alpha,n_\alpha)$:
\begin{align}
	\pm\alpha_2&\longleftrightarrow \pm\left(-\frac{1}{2},\frac{3}{2}\right)\label{WeightTable1}\\
	\pm\alpha_1\pm\alpha_2&\longleftrightarrow \pm\left(\frac{1}{2},\frac{3}{2}\right).\label{WeightTable2}
\end{align}
We will use the description $(m_\alpha,n_\alpha)$ and the corresponding noncompact roots $\alpha$ to refer to the weights of $\frakp_\CC$ interchangeably. The representation $\frakp_\CC$ of $K$ can be decomposed into a direct sum of two 2 dimensional irreducible representations  $V^{\frac{1}{2},\frac{3}{2}}\oplus V^{\frac{1}{2},-\frac{3}{2}}$ of $U(2)$, where the highest weights are labeled as upper indices, with weights
\begin{align*}
	V^{\frac{1}{2},\frac{3}{2}}:& \left(-\frac{1}{2},\frac{3}{2}\right), \left(\frac{1}{2},\frac{3}{2}\right)\\
	V^{\frac{1}{2},-\frac{3}{2}}:& \left(-\frac{1}{2},-\frac{3}{2}\right), \left(\frac{1}{2},-\frac{3}{2}\right).
\end{align*}
The representation $V^{\frac{1}{2},\frac{3}{2}}$ corresponds to the coordinates $z_1$ and $\bar{z}_1$ in the space of matrices given in (\ref{Z1Z2}), and $V^{\frac{1}{2},-\frac{3}{2}}$ corresponds to the coordinates $z_2$ and $\bar{z}_2$ in (\ref{Z1Z2}). If we use gray for noncompact roots and light gray for compact roots, the direct sum decomposition $\frakg_\CC = \frakk_\CC\oplus V^{\frac{1}{2},\frac{3}{2}}\oplus V^{\frac{1}{2},-\frac{3}{2}}$ can be displayed in the root system of $\fraksl(3,\CC)$:
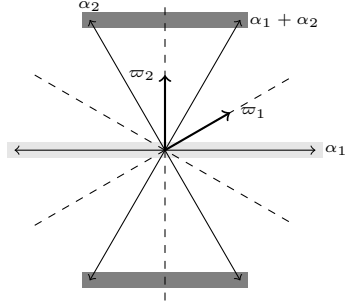
\begin{center}
	    	\begin{tikzpicture}[scale=1]
	\clip (-3,-2) rectangle (3,2);
		\fill[fill=gray!100] (-1.1, 1.632050) rectangle (1.1, 1.832050);
		\fill[fill=gray!100] (-1.1, -1.632050) rectangle (1.1, -1.832050);
		\fill[fill=gray!20] (-2.1,-0.1) rectangle (2.1,0.1);
		\draw[->] (0,0) -- (2,0) node[right] {\tiny $\alpha_1$};
		\draw[->]	(0,0)--(-1,1.732050) node[above] {\tiny $\alpha_2$};
%		\fill[fill=gray!20] (0.05,0.05) --  (1.782050,1.05) -- (0.05,2.05) -- (0.05,0.05);
		\draw[->]	(0,0)--(1,1.732050) node[right] {\tiny $\alpha_1+\alpha_2$};
		\draw[->] (0,0)--(-2,0);
		\draw[->]	(0,0)--(1,-1.732050);
		\draw[->]	(0,0)--(-1,-1.732050);
		\draw[->, thick] (0,0) -- (0.866,0.5)  node[right] {\tiny $\varpi_1$};
		\draw[->, thick] (0,0) -- (0,1)  node[left] {\tiny $\varpi_2$};
		\draw[dashed] (-1.732050,-1) -- (1.732050,1);
		\draw[dashed] (-1.732050,1) -- (1.732050,-1);
		\draw[dashed] (0,-2) -- (0,2);
	\end{tikzpicture}
	\end{center}
There is a weight basis $\{v_{\alpha}\}_{\alpha\in \Delta_{nc}}$ of $\frakp_\CC$, such $v_\alpha$ has weight $\alpha$ when considered as a vector in $\frakg_\CC$ under the adjoint action. By associating these weight vectors to a pair of half integers $(m_\alpha,n_\alpha)$ as above, we can also label each weight vector by writing the pair $(m_\alpha,n_\alpha)$ as lower indices. Therefore the two irreducible $U(2)$ subrepresentations of $\frakp_\CC$ have basis:
\begin{align*}
	V^{\frac{1}{2},\frac{3}{2}}&=\mathrm{span}\{v_{\alpha_2} ,v_{\alpha_1+\alpha_2}\}\\
	V^{\frac{1}{2},-\frac{3}{2}}&=\mathrm{span}\{v_{-\alpha_2} ,v_{-\alpha_1-\alpha_2}\}
\end{align*}
where the vectors $v_\alpha$'s are given by
\begin{align*}
	\begin{matrix}v_{\alpha_2} = v_{-\frac{1}{2}, \frac{3}{2}}=-X_{\alpha_2}& v_{\alpha_1+\alpha_2}=v_{\frac{1}{2}, \frac{3}{2}}=X_{\alpha_1+\alpha_2}\\
	v_{-\alpha_2}=v_{\frac{1}{2}, -\frac{3}{2}}=X_{-\alpha_2}&
	v_{-\alpha_1-\alpha_2}=v_{-\frac{1}{2}, -\frac{3}{2}}=X_{-\alpha_1-\alpha_2}.
	\end{matrix}
\end{align*}

\subsubsection{The maximally noncompact Cartan subalgebra and the restricted roots}
We take the noncompact imaginary root $\alpha_1+\alpha_2\in \Delta_{nc}$ and define the transform
\[
	\ppp _{\alpha_1+\alpha_2} = \Ad\exp(\frac{\pi}{4}(\sigma(v_{\alpha_1+\alpha_2})-v_{\alpha_1+\alpha_2})).
\]
The action by $\ppp_{\alpha_1+\alpha_2}$ sends the complexified maximally compact Cartan subalgebra $\frakt_\CC$  to the maximally noncompact Cartan subalgebra $\fraks_\CC = \frakm_\CC\oplus \fraka_\CC$, where $\fraka$ is the subspace $\CC (X_{\alpha_1+\alpha_2}+X_{-\alpha_1-\alpha_2})$ of $\frakp_\CC$, and $\frakm$ is the centralizer of $\fraka$ in $\frakk$. The analytic subgroup $M\subset G$ with Lie algebra $\frakm$ is
\[
	M=\{e^{-t U_0}e^{3t U_3}=\begin{spmatrix}e^{\ii t}&0&0\\0&e^{-2\ii t}&0\\0&0&e^{\ii t}\end{spmatrix}|t\in\RR\}
\]
The action of the Cayley transform $\qqq_{\alpha_1+\alpha_2} = \ppp_{\alpha_1+\alpha_2}^{-1}$ sends $\fraka_\CC$ to a subspace of $\frakt_\CC$ generated by
\[
	\qqq _{\alpha_1+\alpha_2}(X_{\alpha_1+\alpha_2}+X_{-\alpha_1-\alpha_2}) = H_{\alpha_1}+H_{\alpha_2}
\]
Since $\alpha_i(H_{\alpha_1}+H_{\alpha_2}) = 2$ for $i=1,2$, the two simple roots $\alpha_1, \alpha_2$ have the same restriction to the line $\qqq_{\alpha_1+\alpha_2}\fraka_\CC\subset \frakt_\CC$. We denote this restriction by $\alpha_0$. The set of positive restricted roots $\Delta^+(\frakg,\fraka)$ consists of a character $\alpha_0 = \alpha_1|_{\qqq_{\alpha_1+\alpha_2}\fraka}= \alpha_2|_{\qqq_{\alpha_1+\alpha_2}\fraka}$ with multiplicity 2, and $2\alpha_0 = (\alpha_1+\alpha_2)|_{\qqq_{\alpha_1+\alpha_2}\fraka}$ with multiplicity 1. 
The half sum of positive restricted roots is:
\[
	\rho_0 = \frac{1}{2}\sum_{\alpha\in\Delta^+(\frakg,\fraka) }\alpha = 2\alpha_0
\]
The restricted Weyl group $W(\frakg,\fraka) = N_K(\fraka)/Z_{K}(\fraka)$ is an order 2 group generated by the single element
\[
w_0 = \exp(2 \pi U_3)=\diag(-1,-1,1)
\]
The adjoint action of $\fraka$ on $\frakg$ decomposes the real Lie algebra into restricted root spaces $\frakg_{\pm\alpha_0}$ and $\frakg_{\pm 2\alpha_0}$:
\begin{align*}
	\frakg_{\alpha_0} &= \{\ppp_{\alpha_1+\alpha_2}(z X_{\alpha_1}+\bar{z} X_{\alpha_2})=\frac{1}{\sqrt{2}}\begin{spmatrix}0&z&0\\-\bar{z}&0&\bar{z}\\0&z&0\end{spmatrix}|z\in\CC\}\\
	\frakg_{2\alpha_0} &=\{\ppp_{\alpha_1+\alpha_2}(\ii w X_{\alpha_1+\alpha_2})=-\frac{1}{2}\begin{spmatrix}\ii w&0&-\ii w\\0&0&0\\\ii w&0&-\ii w\end{spmatrix}|w\in\RR\}\\
	\frakg_{-\alpha_0} &= \{\ppp_{\alpha_1+\alpha_2}(\bar{z} X_{-\alpha_1}+z X_{-\alpha_2})=-\frac{1}{\sqrt{2}}\begin{spmatrix}0&z&0\\-\bar{z}&0&-\bar{z}\\0&-z&0\end{spmatrix}|z\in\CC\}\\
	\frakg_{-2\alpha_0} &= \{\ppp_{\alpha_1+\alpha_2}(\ii wX_{-\alpha_1-\alpha_2})=\frac{1}{2}\begin{spmatrix}-\ii w&0&-\ii w\\0&0&0\\\ii w&0&\ii w\end{spmatrix}|w\in\RR\}
\end{align*}
A general element in positive restricted root space $\frakn^+ = \frakg_{\alpha_0}\oplus \frakg_{2\alpha_0}$ can be denoted by
\[
	n_{z,w} = \begin{spmatrix}\ii w&z&-\ii w\\-\bar{z}&0&\bar{z}\\\ii w&z&-\ii w\end{spmatrix}, z\in\CC, w\in\RR
\]
Using the restricted root spaces and the $n_{z,w}$'s, we can compute the Iwasawa decomposition $G=KMA_0N^+$ and the corresponding Iwasawa decomposition $\mathfrak{g}=\mathfrak{k}+\mathfrak{m}+\mathfrak{a}_0+\mathfrak{n}^+$ on the Lie algebra of the $v_\alpha$'s in the complexified Lie algebra $\frakg_\CC$:
\begin{align}
	v_{\alpha_1+\alpha_2}&= -\frac{\ii}{2}(U_0+U_3)+\frac{1}{2}(X_{\alpha_1+\alpha_2}+X_{-\alpha_1-\alpha_2})+\frac{\ii}{2}n_{0,1}\label{IwasawaSU21:1}\\
	v_{-\alpha_1-\alpha_2}&=\frac{\ii}{2}	(U_0+U_3)+\frac{1}{2}(X_{\alpha_1+\alpha_2}+X_{-\alpha_1-\alpha_2})-\frac{\ii}{2}n_{0,1}\label{IwasawaSU21:2}\\
	v_{\alpha_2} &=\ii (U_1-\ii U_2)-\frac{1}{2}n_{1,0} -\frac{\ii}{2} n_{\ii,0}\label{IwasawaSU21:3}\\
	v_{-\alpha_2} &= \ii (U_1+\ii U_2)+\frac{1}{2}n_{1,0}-\frac{\ii}{2} n_{\ii,0}\label{IwasawaSU21:4}
\end{align}
Moreover, the adjoint action of the raising-lowering operators $U_1\pm\ii U_2$ on the $v_\alpha$'s satisfies:
\begin{align*}
	\ad(U_1\pm \ii U_2)v_{\alpha}&= -\ii v_{\alpha_1\pm\alpha}
\end{align*}
if $\alpha_1\pm\alpha$ is a root.

\section{The Principal Series Representations}
Consider the minimal parabolic subgroup $P = MA_0N$, where $A_0$ is the analytic subgroup of $G$ with Lie algebra $\fraka = \RR (X_{\alpha_1+\alpha_2}+X_{-\alpha_1-\alpha_2})$. We introduce the following parameters:
\begin{enumerate}
	\item A character $\fraka\longrightarrow \CC$ sending $X_{\alpha_1+\alpha_2}+X_{-\alpha_1-\alpha_2}$ to a complex number $\lambda$; 
	\item A character $M\longrightarrow\CC$ sending $e^{-t(U_0-3U_3)}$ to $e^{\ii \delta t}$, where $\delta$ is an integer.
\end{enumerate}
Applying the Cayley transform on $\fraks = \frakm\oplus\fraka$,  $\qqq_{\alpha_1+\alpha_2}^{-1}(\fraks)$ is the real Lie algebra $\frakh = \RR\ii (H_{\alpha_1}-H_{\alpha_2})\oplus \RR(H_{\alpha_1}+H_{\alpha_2})$, on which the induction parameter $(\delta,\lambda)$ defines a complex valued character $\chi_{\delta, \lambda}: \frakh\longrightarrow\CC$ such that
\[
	\chi_{\delta, \lambda}(H_{\alpha_1}-H_{\alpha_2}) = \delta, \chi_{\delta, \lambda}(H_{\alpha_1}+H_{\alpha_2}) = \lambda
\]
and therefore the action of $\chi_{\delta, \lambda}$ on $H_{\alpha_i}$'s gives:
\[
	\begin{matrix}\chi_{\delta, \lambda}(H_{\alpha_1})=\frac{\lambda+\delta}{2}, &\chi_{\delta, \lambda}(H_{\alpha_2})=\frac{\lambda-\delta}{2}.\end{matrix}
\]
The character $\chi_{\delta, \lambda}$ lives in the weight lattice if and only if $\lambda\pm\delta\in 2\ZZ$.\\

There are two Casimir elements $\Omega_2$ and $\Omega_3$ of the universal enveloping algebra $U(\frakg_\CC)$, having degree $2$ and $3$ respectively. Denote $\{X_i\}$ to be an indexed basis of $\frakg_\CC$ and $\{\tilde{X}_i\}$ the indexed dual basis under the Killing form. If we take $\pi_{\mathrm{std}}$ to be the standard representation of $\frakg_\CC=\fraksl(3,\CC)$ on a 3 dimensional complex vector space, the quadratic Casimir element can be computed using the formula:
\begin{align*}
	\Omega_2 &= \sum_{i,j}\Tr (\pi_{\mathrm{std}}(X_i)\pi_{\mathrm{std}}(X_j)) \tilde{X}_i\tilde{X}_j\\
	&=\frac{1}{54}(H_{\alpha_1}^2+H_{\alpha_1} H_{\alpha_2}+H_{\alpha_2}^2) + \frac{1}{36}\sum_{\alpha\in \Delta^+}\{X_{\alpha},X_{-\alpha}\}\\
	&=\frac{1}{9}(H_{\alpha_1}^2+H_{\alpha_1} H_{\alpha_2}+H_{\alpha_2}^2+3(H_{\alpha_1} +H_{\alpha_2}))+\frac{1}{18}\sum_{\alpha\in \Delta^+}X_{-\alpha}X_{\alpha}
\end{align*}
where for $X,Y\in U(\frakg_\CC)$, we denote by $\{X,Y\} = X Y+Y X \in U(\frakg_\CC)$. 
Under the Harish-Chandra isomorphism $\gamma': Z(\frakg_\CC)\longrightarrow S(\frakh)^W$, the image of $\Omega_2$ is 
\begin{align*}
	\gamma'(\Omega_2) &= \frac{1}{9}(H_{\alpha_1}^2+H_{\alpha_1} H_{\alpha_2}+H_{\alpha_2}^2-3)
\end{align*}
The cubic Casimir element $\Omega_3$ is:
\begin{align*}
	\Omega_3 &=\sum_{i,j,k}\Tr (\pi_{\mathrm{std}}(X_i)\pi_{\mathrm{std}}(X_j)\pi_{\mathrm{std}}(X_k)) \tilde{X}_i\tilde{X}_j\tilde{X}_k\\
	&=\frac{1}{1944}(3-H_{\alpha _1}+H_{\alpha _2}) (6+2 H_{\alpha _1}+H_{\alpha
   _2}) (H_{\alpha _1}+2 H_{\alpha _2}) \\
   &+\frac{1}{216} (6 X_{-\alpha_1-\alpha_2}X_{\alpha_1+\alpha_2} + 6 X_{-\alpha_2}X_{\alpha_2} - X_{-\alpha_1} (H_1+2H_2) X_{\alpha_1} \\
   &+ X_{-\alpha_2} (2H_1+H_2) X_{\alpha_2}- X_{-\alpha_1-\alpha_2}(H_1-H_2)X_{\alpha_1+\alpha_2} \\
&- 3 X_{-\alpha_2} X_{-\alpha_1} X_{\alpha_1+\alpha_2}
   - 3 X_{-\alpha_1-\alpha_2} X_{\alpha_1}X_{\alpha_2})
\end{align*}
The image of $\Omega_3$ under the Harish-Chandra homomorphism is:
\begin{equation}
	\gamma'(\Omega_3) = -\frac{1}{2^33^5}(H_{\alpha_1}+2H_{\alpha_2}-3)(2H_{\alpha_1}+H_{\alpha_2}+3)(H_{\alpha_1}-H_{\alpha_2}-3)\label{HarishChandraSU21}
\end{equation}
After applying the character $\chi_{\delta, \lambda}$ on $S(\frakh_\CC)^W$, we see that
\begin{align}
	\chi_{\delta, \lambda}(\gamma'(\Omega_2)) &= \frac{1}{36}(3(\lambda^2-4)+\delta^2)\\
	\chi_{\delta, \lambda}(\gamma'(\Omega_3)) &= \frac{1}{2^53^5}(\delta-3)(\delta-3(\lambda-2))(\delta+3(\lambda-2))
\end{align}

We can induce the character $\chi_{\delta, \lambda}$ from the Levi subgroup $L = MA_0$ of the minimal parabolic subgroup $P = MA_0 N$ to get the minimal principal series representation:
\[
	I(\chi_{\delta,\lambda})=\{f:G\rightarrow\CC|f(ge^{-t (U_0-3U_3)}e^{s (X_{\alpha_1+\alpha_2}+X_{-\alpha_1-\alpha_2})}n)=e^{-\ii\delta t-(\lambda+2)s}f(g)\}
\]
The value of the functions in the principal series is determined by their restriction to the maximal compact subgroup $K=U(2)$, and the Lie algebra $\frakg$ acts as differential operators on the left. The functions on $K$ can be expanded into Fourier series with respect to the basis $W^{(j,n)}_{m_1,m_2}$. We can apply the Iwasawa decomposition of the Lie algebra $\frakg = \frakk\oplus\frakm\oplus\fraka\oplus\frakn$ in (\ref{IwasawaSU21:1})-(\ref{IwasawaSU21:4}) and the differential operators $\dd l(\gamma_i)$ and $\dd r(\gamma_i)$ on $K$ introduced in Section {\ref{Irrep}} to write down the action of $\frakg$ on $C^\infty(K)$. Moreover, we can apply the product formula (\ref{ProductClebschGordan}) to express the product of Wigner $D$-functions as linear combinations of Wigner $D$-functions. We will use the machinery developed before to prove the following proposition describing the $\frakg$-action on $I(\chi_{\delta,\lambda})$ explicitly in the next two sections:
\begin{proposition}\label{SU21PrincipalSeriesProp}
	Let $v_\alpha$ be the weight vectors of $\frakp_\CC$ as an $U(2)$ representation, such that $\alpha\in \Delta_{nc}$ as listed in (\ref{WeightTable1})-(\ref{WeightTable2}). If $\alpha\in\Delta^\pm_{nc}$, the action of the weight vectors $v_\alpha$ in $\frakp_\CC$ satisfies:
\begin{align}
	&\dd l(v_\alpha) W^{(j,n)}_{m_1,m_2}  = \nonumber\\
	&\frac{1}{2\sqrt{2j+1}}\sum_{j_0\in\{\pm\frac{1}{2}\}}
	\left(\begin{smallmatrix}j+j_0,m_1+m_\alpha\\J,m_1,\frac{1}{2},m_\alpha \end{smallmatrix}\right)q_{j_0,\pm}\kappa_{j_0,\pm}(j,n,m_1;\lambda)W^{(j+j_0,n\pm \frac{3}{2})}_{m_1+m_\alpha,m_2\pm \frac{1}{2}},\label{SU21PrincipalSeries}
\end{align}
with the coefficients as shown in the following tables:
\begin{center}
	\begin{tabular}{ccc}
		$q_{j_0,\pm}$& - & +\\\hline
		$j_0=-\frac{1}{2}$ & $\sqrt{j-m_2}$ & $\sqrt{j+m_2}$ \\
		$j_0 = \frac{1}{2}$ & $\sqrt{j+m_2+1}$ &  $\sqrt{j-m_2+1}$\\
	\end{tabular}\\
	\begin{tabular}{ccc}
		$\kappa_{j_0,\pm}$& - & +\\\hline
		$j_0=-\frac{1}{2}$ & $-2j-m_2+n+\lambda$ & $2j-m_2+n-\lambda$ \\ 
		$j_0 = \frac{1}{2}$ & $2j-m_2+n+\lambda+2$ &  $2j+m_2-n+\lambda+2$\\
	\end{tabular}
\end{center}
\end{proposition}

\subsection{Embedding of Principal Series in $C^\infty(K)$}
The $(\frakg,K)$ module of $I(\chi_{\delta, \lambda})$ can be embedded into the space $L^2(K)$ as the subspace consisting of $K$-finite functions $f$ satisfying:
\[
	f(ke^{-t(U_0-3U_3)}) = e^{-\ii\delta t}f(k), \text{ for all }t\in\RR.
\]
This requires that $f$ is a finite linear combination of Wigner $D$-functions $W^{(j,n)}_{m_1,m_2}$ satisfying the condition
\[
	e^{\ii (n-3 m_2) t} = e^{-\ii\delta t}, \text{ for all }t\in\RR
\]
which is equivalent to the condition:
\begin{equation}
	-n+3m_2 = \delta
\end{equation}
Therefore as a vector space, the $(\frakg,K)$ module of the principal series $I(\chi_{\delta,\lambda})$ can be embedded into $C_{\delta}(K)\subset C^\infty(K)$ as an algebraic direct sum:
\[
	I(\chi_{\delta,\lambda}) \subset C_{\delta}(K) :=\bigoplus_{-n+3m_2 = \delta}\CC W^{(j,n)}_{m_1,m_2}\subset C^{\infty}(K).
\]
Since $m_2$ always satisfies $-j\leq m_2\leq j$, $n$ can only take half integer values in the interval $-3j-\delta\leq n\leq 3j-\delta$. Denote the set of the pairs $(j,n)$ satisfying this condition as
\begin{align}
	\mathtt{KTypes}(\delta) &= \{(j,n)\in\frac{1}{2}\ZZ\times\frac{1}{2}\ZZ: -3j-\delta\leq n\leq 3j-\delta\}.
\end{align}
The set $\mathtt{KTypes}(\delta)$ parametrizes all the $K$-isotypic components of the principal series $I(\chi_{\delta,\lambda})$:
\[
	I(\chi_{\delta,\lambda}) = \bigoplus_{(j,n)\in \mathtt{KTypes}(\delta) }\tau^{(j,n)},
\]
where $\tau^{(j,n)}$ is a direct sum of copies of irreducible representations of $U(2)$ with highest weight $(j,n)$. The $K$-isotypic subspaces $\tau^{(j,n)}$ can be decomposed into the direct sum:
\begin{align}
	\tau^{(j,n)} = \bigoplus_{\substack{m_1\in\{-j,-j+1,\ldots,j\}\\m_2 \in \mathtt{M}(j,n,\delta)}} \CC W^{(j,n)}_{m_1,m_2}\label{TauSU21}
\end{align}
where the set $\mathtt{M}(j,n,\delta)$ is defined as:
\begin{align}
	\mathtt{M}(j,n,\delta) &= \{m_2\in\{-j,-j+1,\ldots,j\}: m_2 = \frac{n+\delta}{3}\}.
\end{align}
Since for each $(j,n)\in\mathtt{KTypes(\delta)}$ we have $|\mathtt{M}(j,n,\delta)|=1$ or 0, each $K$-type of the $SU(2,1)$ principal series $I(\chi_{\delta, \lambda})$ has multiplicity at most 1. The $K$-types of $I(\chi_{\delta, \lambda})$, with their multiplicities taken into account, can be displayed on the cone in a subset of $(\frac{1}{2}\ZZ)^3$ with coordinates $(j,n,m_2)$ such that $(j,n)\in\mathtt{KTypes}(\delta)$ and $m_2\in\mathtt{M}(j,n,\delta)$.\\

From this embedding of $I(\chi_{\delta, \lambda})$ into $C^\infty(K)$, the action by the Lie algebra elements of $\frakg_\CC$ can be realized as differential operators on $K$. For any $f\in C_\delta^\infty(K)$, we can extend the domain of $f$ to a vector in $I(\chi_{\delta, \lambda})$ by applying the Iwasawa decomposition $G = KMA_0N$ and the transformation rule of the principal series. More precisely, the actions $\dd l(X)$ and $\dd r(X)$ of $X\in\frakg$ on $f\in C_\delta^\infty(K)\subset I(\chi_{\delta, \lambda})$ are given by
\[
	\begin{matrix}(\dd l(X)f)(k) = \frac{\dd}{\dd t}|_{t=0} f(e^{-tX}k),& (\dd r(X)f)(k) = \frac{\dd}{\dd t}|_{t=0} f(ke^{tX})\end{matrix}
\]
We recall that in Section {\ref{Irrep}}, the extension of $\dd r$ to $\frakk_\CC = \frakk\otimes \CC$ is linear. We can follow the same rule and extend $\dd l$ and $\dd r$ further to $\frakg_\CC$ by setting $\dd l (z X) = z \dd l(X)$ and $\dd r (z X) = z \dd r(X)$ for any $z\in\CC$.
Under such extension, the right action $\dd r(n_{z,w})$ by elements $n_{z,w}\in \frakn^+$ in the nilpotent radical sends any $f\in I(\chi_{\delta, \lambda})$ to 0. The element $X_{\alpha_1+\alpha_2}+X_{-\alpha_1-\alpha_2}\in \fraka$ acts by scalar multiplication:
\begin{align}
	\dd  r(X_{\alpha_1+\alpha_2}+X_{-\alpha_1-\alpha_2})W^{(j,n)}_{m_1, m_2} &= -(\lambda+2) W^{(j,n)}_{m_1, m_2}\label{RightActionTorusSU21}.
\end{align}
Combining (\ref{RightActionTorusSU21}), the action of the Lie algebra $\frakk$ on the Wigner $D$-functions (\ref{CompactAction:1})-(\ref{CompactAction:3}), and the Iwasawa decomposition (\ref{IwasawaSU21:1})-(\ref{IwasawaSU21:4}) for the basis vectors $v_{\alpha}$ in $\frakp_\CC$, the formulae of the right action of $v_{\alpha}$  on $L^2(K)$ are given by
\begin{align}
	\dd r(v_{\pm(\alpha_1+\alpha_2)})W^{(j,n)}_{m_1, m_2}&=\left(\mp\frac{1}{2}\dd r(\ii(U_0+U_3))-\frac{1}{2}(\lambda+2)\right)W^{(j,n)}_{m_1, m_2}\nonumber\\
	&=\frac{1}{2}(\mp n\mp m_2-\lambda-2)W^{(j,n)}_{m_1, m_2}\label{right1}\\
	\dd r(v_{\pm\alpha_2})W^{(j,n)}_{m_1, m_2}&=\dd r(\ii(U_1\mp\ii U_2))W^{(j,n)}_{m_1, m_2}= -\sqrt{(j\mp m_2)(j\pm m_2+1)}W^{(j,n)}_{m_1, m_2\pm 1}\label{right2}.
\end{align}
We can express the left action $\dd l(X)$ by any Lie algebra element $X\in \frakg$ in terms of $\dd r$ using the adjoint action of $K$: \[\dd l(X) = \dd r(-\Ad^{-1}(k)X).\]
Then for any $\alpha\in\Delta^\pm_{nc} = \{\pm\alpha_2,\pm\alpha_1\pm\alpha_2\}$, recalling the correspondence (\ref{WeightTable1})-(\ref{WeightTable2}) of $\alpha$ with the pair of integers $(m_\alpha,n_\alpha)$ and the definition of Wigner $D$-functions as matrix coefficients in (\ref{WignerDefinition}), the left action of $v_\alpha$ on the functions in $C^\infty(K)$ can be expressed as a linear combination of right actions by vectors $v_\alpha$ with $\alpha\in\Delta^\pm_{nc}$, having Wigner $D$-functions $- W^{(\frac{1}{2}, n_{\alpha})}_{m_{\beta},m_\alpha}(k^{-1})$ as coefficients: 
\begin{align*}
	\dd l(v_\alpha) &= \dd r(-\Ad(k^{-1})v_{\alpha})= \sum_{\beta\in\left\{\substack{\Delta^+_{nc}\text{ if }\alpha\in\Delta^+_{nc}\\\Delta^-_{nc}\text{ if }\alpha\in\Delta^-_{nc}}\right.}\dd r( - \overline{W^{(\frac{1}{2}, n_{\alpha})}_{m_{\beta},m_\alpha}(k^{-1})}v_{\beta}).
\end{align*}
The same method for $SL(3,\RR)$ has been provided in \cite{MillerButtcane}. By the unitarity of Wigner $D$-function matrices (\ref{Unitarity1}) and (\ref{UnitarityWigner}), we can change the argument from $k^{-1}$ to $k$ and rearrange the upper and lower indices of Wigner $D$-functions:
\begin{align*}
	\dd l(v_\alpha) &= -\sum_{\beta\in\left\{\substack{\Delta^+_{nc}\text{ if }\alpha\in\Delta^+_{nc}\\\Delta^-_{nc}\text{ if }\alpha\in\Delta^-_{nc}}\right.}W^{(\frac{1}{2}, n_{\alpha})}_{m_{\alpha},m_\beta}(k)\dd r(v_{\beta}).
\end{align*}
%Since $\frakp_\CC^+$ and $\frakp_\CC^-$ are irreducible representations, and each $\beta$ in the sum lies either in $\Delta^+_{nc}(\frakk_\CC,\frakt_\CC)$ or $\Delta^-_{nc}(\frakk_\CC,\frakt_\CC)$, depending on whether $\alpha$ is positive or negative.
After listing all the $\beta's$ in $\Delta^\pm_{nc}$, the sum over $\beta$ above has only two terms. For $\alpha\in\Delta^\pm_{nc}$, we have
\begin{align*}
	\dd l(v_{\alpha}) 
	=&-\left(W^{(\frac{1}{2},\pm\frac{3}{2})}_{m_{\alpha},\mp\frac{1}{2}}(k) \dd r(v_{\pm\alpha_2})+W^{(\frac{1}{2},\pm\frac{3}{2})}_{m_{\alpha},\pm\frac{1}{2}}(k) \dd r(v_{\pm(\alpha_1+\alpha_2)})\right) .
\end{align*}
Applying the formulas for the right action (\ref{right1}) and (\ref{right2}), the left action of $v_\alpha$ on $W^{(j,n)}_{m_1,m_2}$ can be written in terms of products of Wigner $D$-functions:
\begin{align*}
	&\dd l(v_{\alpha})W^{(j,n)}_{m_1, m_2}= \\& -\left(-\sqrt{(j\mp m_2)(j\pm m_2+1)}W^{(\frac{1}{2},\pm\frac{3}{2})}_{m_{\alpha},\mp\frac{1}{2}}W^{(j,n)}_{m_1, m_2\pm 1}+\right.\\&
	\left.\frac{1}{2}(\mp n\mp m_2-\lambda- 2)W^{(\frac{1}{2},\pm\frac{3}{2})}_{m_{\alpha},\pm\frac{1}{2}} W^{(j,n)}_{m_1, m_2}\right)
\end{align*}
Recall from (\ref{ProductClebschGordan}) that the product of Wigner $D$-functions is in fact the linear combination of Wigner $D$-functions for the constituents of the tensor product representations, and the coefficients of this linear combination are products of the Clebsch-Gordan coefficients:
\begin{align*}
	W^{(j_1,n_1)}_{m_1,m_2}W^{(j_2,n_2)}_{m_3,m_4}= \sum_{J\in\{j+\frac{1}{2},j-\frac{1}{2}\}}\left(\begin{smallmatrix}J,m_1+m_3\\j_1,m_1,j_2,m_3\end{smallmatrix}\right)\left(\begin{smallmatrix}J,m_2+m_4\\j_1,m_2,j_2,m_4\end{smallmatrix}\right) W^{(J,n_1+n_2)}_{m_1+m_3,m_2+m_4}
\end{align*}
We can thus combine all the matrix coefficients belonging to the same $J$ in the formula for the left action of any $v_\alpha$ with $\alpha\in\Delta^\pm_{nc}$:
%\begin{align*}
%	\dd l(v_{\alpha})W^{\mu}_{\nu_1, \nu_2} &=\sum_{\mu_3}(\sqrt{(j+m_1)(j-m_1+1)} \left(\begin{smallmatrix}\mu_3,\nu_1+\alpha_2\\\mu,\nu_1-\alpha_1,\alpha_1+\alpha_2,\alpha_1+\alpha_2\end{smallmatrix}\right)-\frac{1}{2}\langle{\alpha}_0,\mu-\lambda-\rho_0\rangle \left(\begin{smallmatrix}\mu_3,\nu_1+\alpha_2\\\mu,\nu_1,\alpha_1+\alpha_2,\alpha_2\end{smallmatrix}\right))\\ 
%	&\left(\begin{smallmatrix}\mu_3,\nu_2+w_{\alpha_1}\alpha\\\mu,\nu_2,\alpha_1+\alpha_2,w_{\alpha_1}\alpha\end{smallmatrix}\right) W^{\mu_3}_{\nu_1+\alpha_2,\nu_2+w_{\alpha_1}\alpha}\\
%	\dd l(v_{\alpha})W^{\mu}_{\nu_1, \nu_2} &=
%\end{align*}
\begin{align}
	&\dd l(v_{\alpha})W^{(j,n)}_{m_1, m_2} =\nonumber\\
	&\sum_{j_0\in\{\pm\frac{1}{2}\}}\left(\sqrt{(j\mp m_2)(j\pm m_2+1)} \left(\begin{smallmatrix}j+j_0,m_2\pm\frac{1}{2}\\j,m_2\pm 1,\frac{1}{2},\mp\frac{1}{2}\end{smallmatrix}\right)
	- \frac{1}{2}(\mp n\mp m_2-\lambda- 2) \left(\begin{smallmatrix}j+j_0,m_2\pm\frac{1}{2}\\j,m_2,\frac{1}{2},\pm\frac{1}{2}\end{smallmatrix}\right)\right)\nonumber\\
&\left(\begin{smallmatrix}j+j_0,m_1+m_\alpha\\j,m_1,\frac{1}{2},m_\alpha\end{smallmatrix}\right)W^{(j+j_0,n\pm\frac{3}{2})}_{m_1+m_\alpha,m_2\pm\frac{1}{2}}.\label{SU21BeforeCG}
\end{align}
If $m_2\neq \pm j$, recall from Table \ref{TableCG1} that the table of Clebsch-Gordan coefficients $\begin{spmatrix}j+j_0,m_2+m_0\\j,m_2,\frac{1}{2},m_0\end{spmatrix}$ for $j_0$ and $m_0$ taking the values $\pm\frac{1}{2}$ is
\begin{center}
	\begin{tabular}{ccc}
		$\begin{spmatrix}j+j_0,m_2+m_0\\j,m_2,\frac{1}{2},m_0\end{spmatrix}$& $m_0=-\frac{1}{2}$ & $m_0=+\frac{1}{2}$\\\hline
		$j_0=-\frac{1}{2}$ & $\sqrt{\frac{j+m_2}{2j+1}}$ & $-\sqrt{\frac{j-m_2}{2j+1}}$ \\
		$j_0 = \frac{1}{2}$ & $\sqrt{\frac{j-m_2+1}{2j+1}}$ &  $\sqrt{\frac{j+m_2+1}{2j+1}}$
	\end{tabular}.
\end{center}
Plugging these Clebsch-Gordan coefficients into the formula (\ref{SU21BeforeCG}) for $\dd l(v_\alpha)$, for each $\alpha\in \Delta^\pm_{nc}$, the action of weight vectors $v_\alpha$ in $\frakp_\CC$ satisfies:
\begin{align*}
	&\dd l(v_\alpha) W^{(j,n)}_{m_1,m_2}  = \frac{1}{2\sqrt{2j+1}}\sum_{j_0\in\{\pm\frac{1}{2}\}}
	\left(\begin{smallmatrix}j+j_0,m_1+m_\alpha\\j,m_1,\frac{1}{2},m_\alpha \end{smallmatrix}\right)\times \\&q_{j_0,\pm}(j,m_2)\kappa_{j_0,\pm}(j,n,m_2;\lambda)W^{(j+j_0,n\pm \frac{3}{2})}_{m_1+m_\alpha,m_2\pm \frac{1}{2}},
\end{align*}
which is the formula (\ref{SU21PrincipalSeries}). The expressions of the coefficients $q_{j_0,\pm}$ and $\kappa_{j_0,\pm}$ are shown in the following tables:
\begin{center}
	\begin{tabular}{ccc}
		$q_{j_0,\pm}(j,m_2)$& - & +\\\hline
		$j_0=-\frac{1}{2}$ & $\sqrt{j+m_2}$ & $\sqrt{j-m_2}$ \\
		$j_0 = \frac{1}{2}$ & $\sqrt{j-m_2+1}$ &  $\sqrt{j+m_2+1}$\\
	\end{tabular}\\
	\begin{tabular}{ccc}
		$\kappa_{j_0,\pm}(j,n,m_2;\lambda)$& - & +\\\hline
		$j_0=-\frac{1}{2}$ & $-(2j-m_2+n-\lambda)$ & $2j+m_2-n-\lambda$ \\ 
		$j_0 = \frac{1}{2}$ & $2j+m_2-n+\lambda+2$ &  $2j-m_2+n+\lambda+2$\\
	\end{tabular}
\end{center}
We have thus finished the proof of the Proposition \ref{SU21PrincipalSeriesProp}.

%\begin{align*}
%	C_2 &=\frac{1}{6}(u_{\alpha_1}u_{-\alpha_1}+u_{-\alpha_1}u_{\alpha_1})-\frac{1}{6}(u_{\alpha_2}u_{-\alpha_2}+u_{-\alpha_2}u_{\alpha_2}) \\
%	&+\frac{1}{6}(u_{\alpha_1+\alpha_2}u_{-\alpha_1-\alpha_2}+u_{-\alpha_1-\alpha_2}u_{\alpha_1+\alpha_2})+\frac{1}{9} (H_{\alpha_1}^2+H_{\alpha_1}H_{\alpha_2}+H_{\alpha_2}^2)
%\end{align*}
%\begin{align*}
%	C_2 &=\frac{1}{6}(u_{\alpha_1}u_{-\alpha_1}+u_{-\alpha_1}u_{\alpha_1})-\frac{1}{6}(u_{\alpha_2}u_{-\alpha_2}+u_{-\alpha_2}u_{\alpha_2}) \\
%	&+\frac{1}{6}(u_{\alpha_1+\alpha_2}u_{-\alpha_1-\alpha_2}+u_{-\alpha_1-\alpha_2}u_{\alpha_1+\alpha_2})-\frac{1}{9}(U_0^2+3 U_3^2)
%\end{align*}

\subsection{Decomposition of $I(\chi_{\delta, \lambda})$}\label{decomp}
In this section, we assume the character $\chi_{\delta,\lambda}$ satisfies $\chi_{\delta,\lambda}(H_{\alpha_i}) \in \ZZ\backslash\{0\}$ for $i = 1,2$. In this case, we are assuming $\lambda$ and $\delta$ will satisfy the condition
\[
	\lambda\pm\delta\in2\ZZ\text{ and } |\lambda-\delta| \geq 2.
\]
Under such assumption, would like to discuss the reducibility and compute the full composition series of the principal series $I(\chi_{\delta, \lambda})$ with $\chi_{\delta,\lambda}$ lying in different open Weyl chambers. The set $\mathtt{M}(j,n;\delta)$ consists of at most one element $m_2 = \frac{n+\delta}{3}$, hence the parameter $m_2$ is completely determined by $n$ and $\delta$ in the expression of $\dd l(v_\alpha)$. It should be pointed out that the same result can be obtained from understanding the order of zeros of the intertwining operator $A(\chi_{\delta,\lambda})$ from the next section, at those points where $\lambda\pm\delta\in2\ZZ$.

The formulas for the coefficients $\kappa_{j_0,\pm}$ of the $\frakp_\CC$ action on $I(\chi_{\delta, \lambda})$ are displayed in the following table:
\begin{center}
	\begin{tabular}{ccc}
		$\kappa_{j_0,\pm}$& - & +\\\hline
		$j_0=-\frac{1}{2}$ & $-\frac{2}{3}(3j+n-\frac{\delta}{2})+\lambda$ & $\frac{2}{3}(3j-n+\frac{\delta}{2})-\lambda$ \\
		$j_0 = \frac{1}{2}$ & $\frac{2}{3}(3j-n+\frac{\delta}{2})+(2+\lambda)$ &  $\frac{2}{3}(3j+n-\frac{\delta}{2})+(2+\lambda)$\\
	\end{tabular}
\end{center}
Recall from the description of the representations of $U(2)$ that $j+n\in \ZZ$ and $j\pm m_2 = j\pm
\frac{n+\delta}{3}\in\ZZ$ in Section {\ref{Irrep}}, let $(k,l)$ be the unique pair of integers  such that $n =- \delta+ \frac{3}{2}l, j = \frac{k}{2}$. They live in the following cone $\mathtt{LatticeCond}$ of the lattice $\ZZ^2$:
\[(k,l)\in \mathtt{LatticeCond} = \{(k,l)\in \ZZ_{\geq 0}\times \ZZ| -k\leq l\leq k\text{ and }k\equiv l\text{ }\mathrm{mod}\text{ } 2\}.\] The coefficients $\kappa_{j_0,\pm}$ can thus be expressed in terms of $k,l,\lambda,\delta$:
\begin{center}
	\begin{tabular}{ccc}
		$\kappa_{j_0,\pm}$& - & +\\\hline
		$j_0=-\frac{1}{2}$ & $-(k+l-\lambda-\delta)$ & $k-l-\lambda+\delta$ \\
		$j_0 = \frac{1}{2}$ & $k-l+\lambda+\delta+2$ &  $k+l+\lambda-\delta+2$\\
	\end{tabular}.
\end{center}

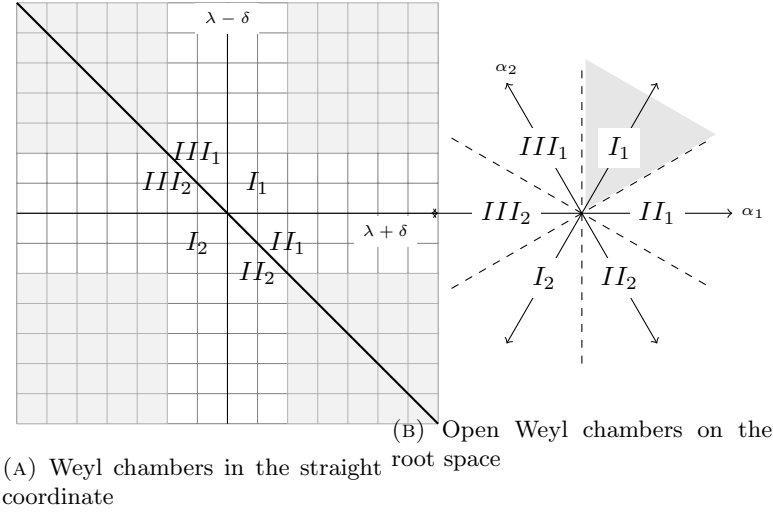
\begin{figure}[h]
  \centering
  \begin{subfigure}{0.4\textwidth}
    \centering
	\begin{tikzpicture}[scale=0.4]
	\clip (-7.5,-7.5) rectangle (7.5,7.5);
	\draw[step=1,gray,very thin] (-7,-7) grid (7,7);
	
	\fill[fill=gray!20,opacity=0.5] (2,2) rectangle (7,7);
	\fill[fill=gray!20,opacity=0.5] (-2,2) rectangle (-7,7);
	\fill[fill=gray!20,opacity=0.5] (2,-2) rectangle (7,-7);
	\fill[fill=gray!20,opacity=0.5] (-2,-2) rectangle (-7,-7);
	\draw[thick] (7,-7) --(-7,7);
	
%	\draw[<-] (4.5,5.5) node[above] {\tiny $V_{\mathrm{fin}}$} -- (5.5,4.5);
%	\draw[<-] (5.5,4.5) node[right] {\tiny $V_{\mathrm{comp}+}$} -- (4.5,3.5);
%	\draw[<-] (3.5,4.5) -- (4.5,3.5) node[below] {\tiny $V_{\mathbb{H}}$};
%	\draw[<-] (4.5,5.5) -- (3.5,4.5) node[left] {\tiny $V_{\mathrm{comp}-}$};
	
%	\draw[->] (-4.5,-5.5) node[below] {\tiny $V_{\mathrm{fin}}$} -- (-5.5,-4.5);
%	\draw[->] (-5.5,-4.5) node[left] {\tiny $V_{\mathrm{comp}+}$} -- (-4.5,-3.5);
%	\draw[->] (-3.5,-4.5) -- (-4.5,-3.5) node[above] {\tiny $V_{\mathbb{H}}$};
%	\draw[->] (-4.5,-5.5) -- (-3.5,-4.5) node[right] {\tiny $V_{\mathrm{comp}-}$};
	
%	\draw[->] (6,-5) node[above,right] {\tiny $V_{\mathrm{disc}+}$} -- +(-1,1) node[above,right] {\tiny $V_{\mathbb{H}}$};
%	\draw[->] (5,-4) -- +(-1,1) node[above,right] {\tiny $V_{\mathrm{comp}+}$};
	
%	\draw[->] (5,-6) node[below,left] {\tiny $V_{\mathrm{comp}+}$} -- +(-1,1) node[below,left] {\tiny $V_{\mathbb{H}}$};
%	\draw[->] (4,-5) -- +(-1,1) node[below] {\tiny $V_{\mathrm{disc}+}$};
	
	\draw (1,1) node { $I_1$};
	\draw (-1,2) node { $III_1$};
	\draw (-2,1) node { $III_2$};
	\draw (-1,-1) node {$I_2$};
	\draw (2,-1) node { $II_1$};
	\draw (1,-2) node { $II_2$};
	
	\draw[->] (-7,0) --  (7,0) node[left=0.7cm, below = 1pt, fill=white] {\tiny $\lambda+\delta$};
	\draw[->] (0,-7) -- (0,7) node[below, fill=white] {\tiny $\lambda-\delta$};
	\end{tikzpicture}
	\caption{Weyl chambers in the straight coordinate}
	\label{paramquad}
  \end{subfigure}
    \begin{subfigure}{0.4\textwidth}
    \centering
    	\begin{tikzpicture}
	\clip (-2.5,-2.5) rectangle (2.5,2.5);
		\draw[->] (0,0) -- (2,0) node[right] {\tiny $\alpha_1$};
		\draw[->]	(0,0)--(-1,1.732050) node[above] {\tiny $\alpha_2$};
		\fill[fill=gray!20] (0.05,0.05) --  (1.782050,1.05) -- (0.05,2.05) -- (0.05,0.05);
		\draw[->]	(0,0)--(1,1.732050);
		\draw[->] (0,0)--(-2,0);
		\draw[->]	(0,0)--(1,-1.732050);
		\draw[->]	(0,0)--(-1,-1.732050);
		\draw[dashed] (-1.732050,-1) -- (1.732050,1);
		\draw[dashed] (-1.732050,1) -- (1.732050,-1);
		\draw[dashed] (0,-2) -- (0,2);
		\draw (0.5,0.866) node[fill=white] {$I_1$};
		\draw (-0.5,0.866) node[fill=white] {$III_1$};
		\draw (-1,0) node[fill=white] {$III_2$};
		\draw (1,0) node[fill=white] {$II_1$};
		\draw (-0.5,-0.866) node[fill=white] {$I_2$};
		\draw (0.5,-0.866) node[fill=white] {$II_2$};
	\end{tikzpicture}
	\caption{Open Weyl chambers on the root space}
  \end{subfigure}
  \caption{Weyl chambers in different coordinates}
  \label{weylchamber}
\end{figure}

%$\lambda-\delta\geq 0, \lambda+\delta\geq 0,\lambda-\delta\leq -2, \lambda+\delta \leq -2$
Based on the signs of $\frac{\lambda\pm\delta}{2}$ and $\lambda$, which are the values of $\chi_{\delta,\lambda}$ on coroots, the dual of Cartan subalgebra $\frakh_\CC^*$ is divided into 6 Weyl chambers as shown in the Figure \ref{weylchamber}. The Weyl group $W$ acts on the characters $\chi_{\delta, \lambda}$, sending it to different Weyl chambers. The action by the simple reflections $w_{\alpha_1},w_{\alpha_2}$ on the pair of parameters $(\delta,\lambda)$ and the corresponding $\chi_{\delta, \lambda}$ is:
\begin{align*}
	w_{\alpha_1}\chi_{\delta,\lambda} &= \chi_{-\frac{3\lambda+\delta}{2}, \frac{\lambda-\delta}{2}}\\
	w_{\alpha_2}\chi_{\delta,\lambda} &= \chi_{\frac{3\lambda-\delta}{2}, \frac{\lambda+\delta}{2}}.
\end{align*}

\subsubsection{Modules}
There are 6 families of irreducible $(\frakg,K)$-modules for the group $SU(2,1)$ depending on the parameters $(\delta,\lambda)$. They can be decomposed into direct sums of $\tau^{(j,n)}$'s as defined in (\ref{TauSU21}).
%\begin{enumerate}
%	\item Holomorphic discrete series $V_{\mathrm{disc}+}$ and antiholomorphic discrete series $V_{\mathrm{disc}-}$:
%	\begin{align*}
%		V_{\mathrm{disc}\pm}(\chi_{\delta, \lambda}) = \bigoplus_{\substack{(k,l)\in \mathtt{LatticeCond}\\k\mp l<\lambda\pm\delta}}\tau^{(\frac{k}{2},\frac{\pm 3\lambda+\delta}{2}+\frac{3l}{2})}
%	\end{align*}
%	\item Quaternionic discrete series $V_{\mathbb{H}}$:
%	\begin{align*}
%		V_{\mathbb{H}}(\chi_{\delta, \lambda})  &= \bigoplus_{\substack{(k,l)\in \mathtt{LatticeCond}\\k-l\leq\lambda-\delta\\k+l\geq \lambda+\delta}} \tau^{(\frac{k}{2}, -\delta+\frac{3l}{2})}
%	\end{align*}
%	\item Finite dimensional representation $V_{\mathrm{fin}}$:
%	\begin{align*}
%		V_{\mathrm{fin}}(\chi_{\delta, \lambda})  &= \bigoplus_{\substack{(k,l)\in \mathtt{LatticeCond}\\k+l<\lambda+\delta\\k-l<\lambda-\delta}} \tau^{(\frac{k}{2},-\delta + \frac{3l}{2})}
%	\end{align*}
%	\item Two other modules $Q_{\pm}$:
%	\begin{align*}
%		Q_{\pm}(\chi_{\delta, \lambda})  &= \bigoplus_{\substack{(k,l)\in \mathtt{LatticeCond}\\  k\mp l < \lambda\mp\delta}} \tau^{(\frac{k}{2},-\delta+\frac{3l}{2})}
%	\end{align*}
%\end{enumerate}
We are going to display these modules in diagrams of lattice points and shaded regions in $(k,l)$-coordinates. In these diagrams, the lattice points stand for $K$-types $\tau^{(j,n)}$ represented in $(k,l)$ coordinates. The horizontal and vertical axes stand for $k$ and $l$, respectively.  The dashed arrows stand for a possible action that maps one $K$-type to another by the Lie algebra action. The irreducible subquotients are depicted by regions of different shades of gray, the darkest gray is for the finite dimensional representation $V_\mathrm{fin}$ or the holomorphic/antiholomorphic discrete series $V_{\mathrm{disc}\pm}$, the medium gray is for $Q_\pm$ and the lightest gray is for the \emph{quaternionic discrete series} $V_{\mathbb{H}}$. In these pictures, the lowest $K$-types are labeled by $(j,n)$ instead of $(k,l)$.
\subsubsection{The Weyl chamber $I_1$}
The character $\chi_{\delta,\lambda}$ in the Weyl chamber $I_1$ satisfies
 \[\lambda-\delta\geq 2, \lambda+\delta\geq 2.\]
There exists $(\frakg,K)$-submodules of $I(\chi_{\delta, \lambda})$ generated by finitely many $K$-types:
\begin{align*}
	V_{\mathbb{H}}(\chi_{\delta, \lambda}) &= U(\frakg) \tau^{(\frac{\lambda}{2}, \frac{\delta}{2})}\\
	V_{1}(\chi_{\delta, \lambda}) &= U(\frakg)\tau^{(\frac{\lambda+\delta}{4}, \frac{3\lambda-\delta}{4})} + U(\frakg) \tau^{(\frac{\lambda-\delta}{4}, \frac{-3\lambda-\delta}{4})}
\end{align*}
that form a composition series of $I(\chi_{\delta, \lambda})$:
\[
	V_{\mathbb{H}}\xhookrightarrow[V_1/V_{\mathbb{H}} = Q_{-}\oplus Q_{+}]{\iota_2} V_1\xhookrightarrow[V_0/V_1 = V_{\mathrm{fin}}]{\iota_1} V_0= I(\chi_{\delta, \lambda}).
\]
The quaternionic discrete series $V_{\mathbb{H}}$, a finite dimensional representation $V_{\mathrm{fin}}$ and the $Q_{\pm}$'s  are irreducible $(\frakg, K)$ modules, which decompose into a direct sum of $K$-isotypic spaces:
\begin{align*}
	V_{\mathbb{H}}(\chi_{\delta, \lambda}) &= \bigoplus_{\substack{(k,l)\in \mathtt{LatticeCond}\\k-l\geq\lambda-\delta\\k+l\geq\lambda+\delta}} \tau^{(\frac{k}{2},-\delta+\frac{3l}{2})}\\
	Q_{\pm}(\chi_{\delta, \lambda}) &= \bigoplus_{\substack{(k,l)\in \mathtt{LatticeCond}\\k\mp l <\lambda\mp\delta\\k\pm l \geq\lambda\pm\delta}{2}} \tau^{(\frac{k}{2},-\delta+\frac{3l}{2})}\\
	V_{\mathrm{fin}} (\chi_{\delta, \lambda})&= \bigoplus_{\substack{(k,l)\in \mathtt{LatticeCond}\\k+l<\lambda+\delta\\k-l<\lambda-\delta }} \tau^{(\frac{k}{2},-\delta + \frac{3l}{2})}
\end{align*}
%For each $\mu$ such that $(j(\mu),n(\mu)) = (j,n)$, $\mu$ can be written as the linear combination $\mu = (j+\frac{n}{3})\alpha_1+\frac{2}{3}n\alpha_2$.
An example when $(\delta,\lambda) = (0,4)$ is displayed in the figure below. 
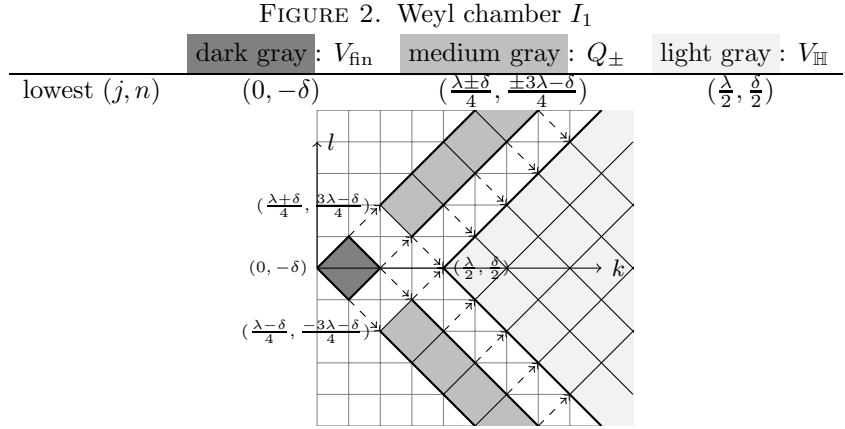
\begin{figure}[h]
\centering
 \caption{Weyl chamber $I_1$}
\begin{tabular}{cccc}&\colorbox{gray}{dark gray}: $V_{\mathrm{fin}}$&\colorbox{gray!50}{medium gray}: $Q_\pm$& \colorbox{gray!10}{light gray}: $V_\mathbb{H}$\\\hline
lowest $(j,n)$ & $(0,-\delta)$ & $(\frac{\lambda\pm\delta}{4},\frac{\pm 3\lambda-\delta}{4})$ & $(\frac{\lambda}{2},\frac{\delta}{2})$
\end{tabular}
\begin{tikzpicture}[scale=0.42]
	\clip (-3.1,-5) rectangle (10,5);
	\draw[step=1,gray,very thin] (0,-5) grid (10,5);

	\draw[fill=gray!10] (4,0) -- (14,10) -- (14,-10);
	\draw[fill=gray!100] (0,0)--(1,1)--(2,0)--(1,-1);
	\draw[fill=gray!50] (2,2)--(5,5)--(8,6)--(3,1);
	\draw[fill=gray!50] (2,-2)--(5,-5)--(8,-6)--(3,-1);
	
	\draw[thick] (2,2) -- (5,5);
	\draw[thick] (2,-2) -- (5,-5);
		\foreach \x in {2,3,4,5,6}
			\draw (\x,\x) -- (\x+1,\x-1);
		\foreach \y in {2,3,4,5,6}
			\draw (\y,-\y) -- (\y+1,-\y+1);
		\foreach \x in {3,4,...,7}
			\draw[dashed,->>] (\x,\x-2) -- (\x+1,\x-3);
		\foreach \y in {3,4,...,7}
			\draw[dashed,->>] (\y,-\y+2) -- (\y+1,-\y+3);
	\foreach \x in {1,2}
		\draw[dashed,->>] (\x,2-\x) -- (\x+1,3-\x);
	\foreach \x in {1,2}
		\draw[dashed,->>] (\x,-2+\x) -- (\x+1,-3+\x);
	\draw[thick] (0,0) -- (1,1)--(2,0)--(1,-1)--(0,0);
	\draw[thick] (3,1) -- (8,6);
	\draw[thick] (3,-1) -- (8,-6);
	\draw[thick] (4,0) -- (9,5);
	\draw[thick] (4,0) -- (9,-5);
		\foreach \x in {5,6,7,8,9}
			\draw (\x-1,5-\x) -- (\x+4,10-\x);
		\foreach \y in {5,6,7,8,9}
			\draw (\y-1,\y-5) -- (\y+4,\y-10);

	\draw (0,0) node [left] {\tiny $(0,-\delta)$};
	\draw (2,2) node [left] {\tiny $(\frac{\lambda+\delta}{4},\frac{3\lambda-\delta}{4})$};
	\draw (2,-2) node [left] {\tiny $(\frac{\lambda-\delta}{4},\frac{-3\lambda-\delta}{4})$};
	\draw (4,0) node [right] {\tiny $(\frac{\lambda}{2},\frac{\delta}{2})$};
	
	\draw[right=1pt] (9,0) node {\small $k$};
	\draw[right=1pt] (0,4) node {\small $l$};
	\draw[->] (0,0) -- (0,4);
	\draw[->] (0,0) -- (9,0);
\end{tikzpicture}
\end{figure}

\subsubsection{The Weyl chamber $II_1$}
The character $\chi_{\delta,\lambda}$ in the Weyl chamber $II_1$ satisfies \[\lambda > 0, \lambda-\delta \leq -2.\]
The two $(\frakg,K)$-submodules of $I(\chi_{\delta, \lambda})$ in the composition series are:
\begin{align*}
	V_{\mathbb{H}}(w_{\alpha_2}\chi_{\delta, \lambda}) &= U(\frakg) \tau^{(\frac{\lambda+\delta}{4}, \frac{3\lambda-\delta}{4})}\\
	V_{\mathrm{disc}-} (w_{\alpha_2}\chi_{\delta, \lambda})&= U(\frakg)\tau^{(0, -\delta)},
%    W_{2} &= U(\frakg) W^{(\frac{-\lambda-\delta}{4}, \frac{-3\lambda+\delta}{4})}_{*,*}\\
\end{align*}
where $V_{\mathbb{H}}$ is the quaternionic discrete series, and $V_{\mathrm{disc}-}$ is the antiholomorphic discrete series. Quotienting out the direct sum of these two modules from $I(\chi_{\delta, \lambda})$, we can get 
\[
	V_{\mathbb{H}}\oplus V_{\mathrm{disc}-} \xhookrightarrow[V_0/(V_{\mathbb{H}}\oplus V_{\mathrm{disc}-}) = Q_-]{\iota} V_0=I(\chi_{\delta, \lambda}).
\]
The spaces $V_{\mathrm{disc}-}, V_{\mathbb{H}}$ and $Q_-$ are irreducible $(\frakg, K)$ modules, which are the direct sum of $K$-types:
\begin{align*}
	V_{\mathbb{H}}(w_{\alpha_2}\chi_{\delta, \lambda}) &= \bigoplus_{\substack{(k,l)\in \mathtt{LatticeCond}\\k+l\geq \lambda+\delta}} \tau^{(\frac{k}{2}, -\delta+\frac{3l}{2})}\\
	Q_-(w_{\alpha_2}\chi_{\delta, \lambda}) &= \bigoplus_{\substack{(k,l)\in \mathtt{LatticeCond}\\ k+l < \lambda+\delta\\k+l \geq -\lambda+\delta}} \tau^{(\frac{k}{2},-\delta+\frac{3l}{2})}\\
	V_{\mathrm{disc}-} (w_{\alpha_2}\chi_{\delta, \lambda})&= \bigoplus_{\substack{(k,l)\in \mathtt{LatticeCond}\\ k+l<-\lambda+\delta}} \tau^{(\frac{k}{2},-\delta + \frac{3l}{2})}
\end{align*}
For $(\delta,\lambda) = (6,2)$, the regions representing the $K$-types of these modules are shown in the picture below:
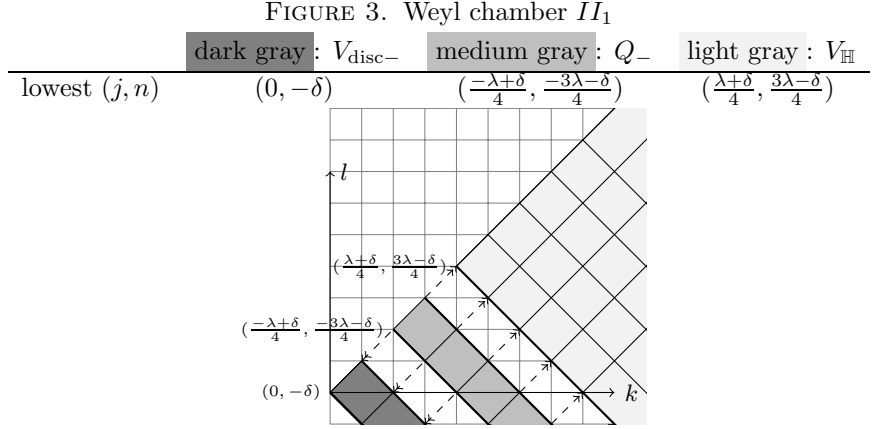
\begin{figure}[h]
\centering
 \caption{Weyl chamber $II_1$}
\begin{tabular}{cccc}
&\colorbox{gray}{dark gray}: $V_{\mathrm{disc}-}$ & \colorbox{gray!50}{medium gray}: $Q_-$& \colorbox{gray!10}{light gray}: $V_\mathbb{H}$\\\hline
lowest $(j,n)$ & $(0,-\delta)$&$(\frac{-\lambda+\delta}{4},\frac{-3\lambda-\delta}{4})$&$(\frac{\lambda+\delta}{4},\frac{3\lambda-\delta}{4})$
\end{tabular}
\begin{tikzpicture}[scale=0.42]
	\clip (-3.1,-5) rectangle (10,5);
	\draw[step=1,gray,very thin] (0,-5) grid (10,5);
	\draw[fill=gray!10] (4,0) -- (14,10) -- (14,-10);
	\draw[fill=gray!50] (2,-2)--(5,-5)--(8,-6)--(3,-1);
	\draw[fill=gray!100] (0,-4)--(1,-3)--(6,-8)--(5,-9);
		\foreach \x in {3,4,...,9}
			\draw[dashed,->>] (\x,-\x+6-4) -- +(1,1);
		\foreach \y in {1,2,...,9}
			\draw[dashed,<<-] (\y,-\y+2-4) -- (\y+1,-\y+3-4);
	\draw[thick] (0,0-4) -- +(5,-5);
	\draw[thick] (1,1-4) -- +(6,-6);
	\draw[thick] (0,0-4) -- +(5,-5);
	\draw[thick] (2,2-4) -- +(7,-7);
	\draw[thick] (3,3-4) -- +(8,-8);
	\draw[thick] (4,4-4) -- +(9,-9);
	\foreach \x in {4,5,...,8}
		\draw (\x, -\x+8-4) -- +(9,9);
	\foreach \x in {5,6,...,8}
		\draw (\x, \x-4) -- +(9,-9);
	\foreach \x in {0,1,...,8}
			\draw (\x,-\x-4) -- +(1,1);
	\foreach \x in {2,3,...,8}
			\draw (\x,-\x+4-4) -- +(1,1);

	\draw (0,-4) node [left] {\tiny $(0,-\delta)$};
	\draw (2,-2) node [left] {\tiny $(\frac{-\lambda+\delta}{4},\frac{-3\lambda-\delta}{4})$};
	\draw (4,0) node [left] {\tiny $(\frac{\lambda+\delta}{4},\frac{3\lambda-\delta}{4})$};
	
	\draw[right=1pt] (9,-4) node {\small $k$};
	\draw[right=1pt] (0,3) node {\small $l$};
\draw[->] (0,0-4) -- (0,7-4);
\draw[->] (0,0-4) -- (9,0-4);
\end{tikzpicture}
\end{figure}

\subsubsection{The Weyl chamber $II_2$}
The character $\chi_{\delta,\lambda}$ lying in the Weyl chamber $II_2$ satisfies the inequality: \[\lambda < 0, \lambda+\delta \geq 2.\]There exists a $(\frakg,K)$-submodule
\begin{align*}
%V_{\mathbb{H}} &= U(\frakg) W^{(\frac{-\lambda-\delta}{4}, \frac{-3\lambda+\delta}{4})}_{*,*}\\
%	V_{2} &= U(\frakg)W^{(0, \delta)}_{*,*} \\
    Q_-(w_{\alpha_1}w_{\alpha_2}\chi_{\delta, \lambda}) &= U(\frakg) \tau^{(\frac{\lambda+\delta}{4}, \frac{3\lambda-\delta}{4})}
\end{align*}
of $I(\chi_{\delta, \lambda})$ that forms a composition series of $I(\chi_{\delta, \lambda})$:
\[
	Q_- \xhookrightarrow[V_0/Q_- = V_{\mathbb{H}}\oplus V_{\mathrm{disc}-}]{\iota} V_0=I(\chi_{\delta,\lambda}).
\]

The spaces $V_{\mathrm{disc}-}, V_{\mathbb{H}}$ and $Q_-$ are irreducible $(\frakg, K)$ modules, which are direct sums of $K$-types:
\begin{align*}
	V_{\mathbb{H}} (w_{\alpha_1}w_{\alpha_2}\chi_{\delta, \lambda})&= \bigoplus_{\substack{(k,l)\in \mathtt{LatticeCond}\\k+l\geq -\lambda+\delta}} \tau^{(\frac{k}{2}, -\delta+\frac{3l}{2})}\\
	Q_- (w_{\alpha_1}w_{\alpha_2}\chi_{\delta, \lambda})&= \bigoplus_{\substack{(k,l)\in \mathtt{LatticeCond}\\ k+l < -\lambda+\delta\\k+l \geq \lambda+\delta}} \tau^{(\frac{k}{2},-\delta+\frac{3l}{2})}\\
	V_{\mathrm{disc}-}(w_{\alpha_1}w_{\alpha_2}\chi_{\delta, \lambda}) &= \bigoplus_{\substack{(k,l)\in \mathtt{LatticeCond}\\ k+l<\lambda+\delta}} \tau^{(\frac{k}{2},-\delta + \frac{3l}{2})}
\end{align*}
For $(\delta,\lambda) = (6,-2)$, the regions representing the $K$-types of these modules are shown in the picture below:
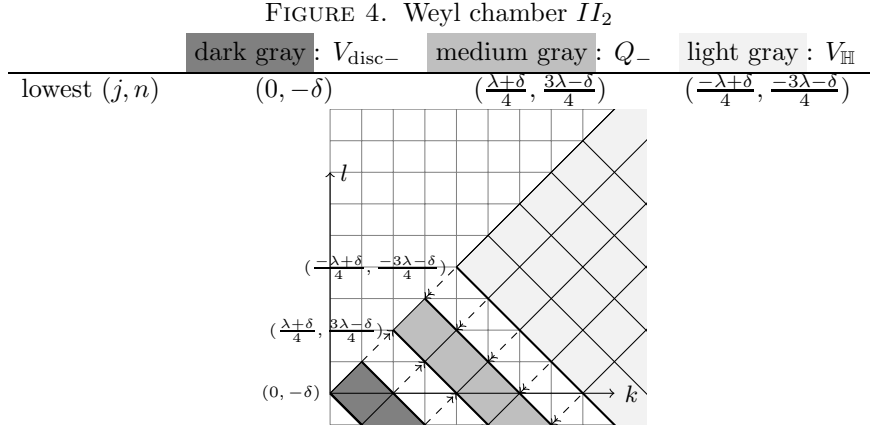
\begin{figure}[h]
\centering
 \caption{Weyl chamber $II_2$}
 \begin{tabular}{cccc}
 &\colorbox{gray}{dark gray}: $V_{\mathrm{disc}-}$&\colorbox{gray!50}{medium gray}: $Q_-$&\colorbox{gray!10}{light gray}: $V_\mathbb{H}$\\\hline
 lowest $(j,n)$ & $(0,-\delta)$ & $(\frac{\lambda+\delta}{4},\frac{3\lambda-\delta}{4})$ & $(\frac{-\lambda+\delta}{4},\frac{-3\lambda-\delta}{4})$
 \end{tabular}
\begin{tikzpicture}[scale=0.42]
	\clip (-3.1,-5) rectangle (10,5);
	\draw[step=1,gray,very thin] (0,-5) grid (10,5);
	\draw[fill=gray!10] (4,0) -- (14,10) -- (14,-10);
	\draw[fill=gray!50] (2,-2)--(5,-5)--(8,-6)--(3,-1);
	\draw[fill=gray!100] (0,-4)--(1,-3)--(6,-8)--(5,-9);
	\foreach \x in {3,4,...,9}
			\draw[dashed,<<-] (\x,-\x+6-4) -- +(1,1);
		\foreach \y in {1,2,...,9}
			\draw[dashed,->>] (\y,-\y+2-4) -- (\y+1,-\y+3-4);
	\draw[thick] (0,0-4) -- +(5,-5);
	\draw[thick] (1,1-4) -- +(6,-6);
	\draw[thick] (0,0-4) -- +(5,-5);
	\draw[thick] (2,2-4) -- +(7,-7);
	\draw[thick] (3,3-4) -- +(8,-8);
	\draw[thick] (4,4-4) -- +(9,-9);
	\foreach \x in {4,5,...,8}
		\draw (\x, -\x+8-4) -- +(9,9);
	\foreach \x in {5,6,...,8}
		\draw (\x, \x-4) -- +(9,-9);
	\foreach \x in {0,1,...,8}
			\draw (\x,-\x-4) -- +(1,1);
	\foreach \x in {2,3,...,8}
			\draw (\x,-\x+4-4) -- +(1,1);

\draw (0,-4) node [left] {\tiny $(0,-\delta)$};
	\draw (4,0) node [left] {\tiny $(\frac{-\lambda+\delta}{4},\frac{-3\lambda-\delta}{4})$};
	\draw (2,-2) node [left] {\tiny $(\frac{\lambda+\delta}{4},\frac{3\lambda-\delta}{4})$};
		
	\draw[right=1pt] (9,-4) node {\small $k$};
	\draw[right=1pt] (0,3) node {\small $l$};
\draw[->] (0,0-4) -- (0,7-4);
\draw[->] (0,0-4) -- (9,0-4);
\end{tikzpicture}
\end{figure}

\subsubsection{The Weyl chamber $I_2$}
When the character $\chi_{\delta,\lambda}$ lies in the Weyl chamber $I_2$,\[\lambda-\delta\leq -2, \lambda+\delta \leq -2.\]There exists $(\frakg,K)$-submodules
\begin{align*}
%	V_{\mathbb{H}} &= U(\frakg) W^{(\frac{\lambda}{2}, -\frac{\delta}{2})}_{*,*}\\
	V_{1}(w_{\alpha_1}w_{\alpha_2}w_{\alpha_1}\chi_{\delta, \lambda}) &= U(\frakg)\tau^{(\frac{-\lambda+\delta}{4}, \frac{-3\lambda-\delta}{4})}+ U(\frakg) \tau^{(\frac{-\lambda-\delta}{4}, \frac{3\lambda-\delta}{4})}\\
    V_{\mathrm{fin}}(w_{\alpha_1}w_{\alpha_2}w_{\alpha_1}\chi_{\delta, \lambda}) &=U(\frakg)\tau^{(0,-\delta)}
\end{align*}
of $I(\chi_{\delta, \lambda})$ that forms a composition series of $I(\chi_{\delta, \lambda})$:
\[
	V_{\mathrm{fin}}\xhookrightarrow[V_1/V_{\mathrm{fin}} = Q_{-}\oplus Q_{+}]{\iota_2} V_1\xhookrightarrow[V_0/V_1 = V_{\mathbb{H}}]{\iota_1} V_0= I(\chi_{\delta,\lambda}).
\]

The spaces $V_\mathrm{fin}, Q_{\pm}$ and $V_{\mathbb{H}}$ are irreducible $(\frakg, K)$ modules, which are direct sums of $K$-types:
\begin{align*}
	V_{\mathbb{H}}(w_{\alpha_1}w_{\alpha_2}w_{\alpha_1}\chi_{\delta, \lambda}) &= \bigoplus_{\substack{(k,l)\in \mathtt{LatticeCond}\\k-l\geq-\lambda+\delta\\k+l\geq-\lambda-\delta}} \tau^{(\frac{k}{2}, -\delta+\frac{3l}{2})}\\
	Q_{\pm} (w_{\alpha_1}w_{\alpha_2}w_{\alpha_1}\chi_{\delta, \lambda})&= \bigoplus_{\substack{(k,l)\in \mathtt{LatticeCond}\\k\mp l <-\lambda\pm\delta \\ k\pm l\geq -\lambda\mp\delta}} \tau^{(\frac{k}{2},-\delta+\frac{3l}{2})}\\
	V_{\mathrm{fin}} (w_{\alpha_1}w_{\alpha_2}w_{\alpha_1}\chi_{\delta, \lambda})&= \bigoplus_{\substack{(k,l)\in \mathtt{LatticeCond}\\k-l<-\lambda+\delta\\k+l<-\lambda-\delta}} \tau^{(\frac{k}{2},-\delta + \frac{3l}{2})}.
\end{align*}
For $(\delta,\lambda) = (0,-4)$, the regions representing the $K$-types of these modules are shown in the picture below:
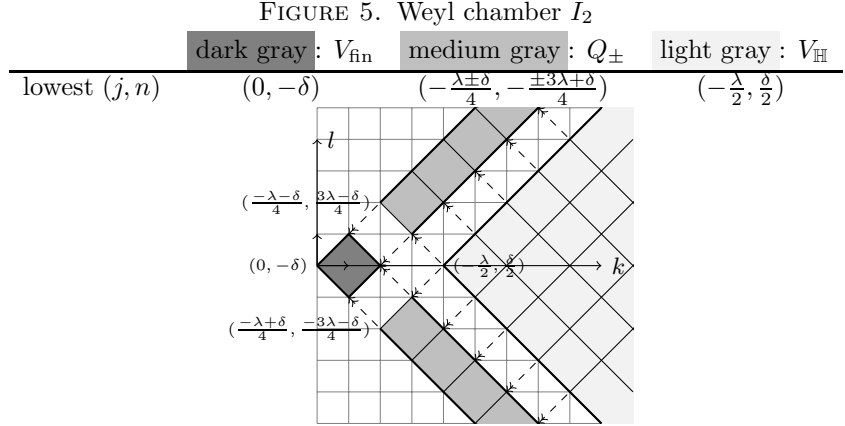
\begin{figure}[h]
\centering
 \caption{Weyl chamber $I_2$}
 \begin{tabular}{cccc}&\colorbox{gray}{dark gray}: $V_{\mathrm{fin}}$&\colorbox{gray!50}{medium gray}: $Q_\pm$& \colorbox{gray!10}{light gray}: $V_\mathbb{H}$\\\hline
lowest $(j,n)$ & $(0,-\delta)$ & $(-\frac{\lambda\pm\delta}{4},-\frac{\pm 3\lambda+\delta}{4})$ & $(-\frac{\lambda}{2},\frac{\delta}{2})$
\end{tabular}
\begin{tikzpicture}
[scale=0.42]
	\clip (-3.1,-5) rectangle (10,5);
	\draw[step=1,gray,very thin] (0,-5) grid (10,5);

	\draw[fill=gray!10] (4,0) -- (14,10) -- (14,-10);
	\draw[fill=gray!100] (0,0)--(1,1)--(2,0)--(1,-1);
	\draw[fill=gray!50] (2,2)--(5,5)--(8,6)--(3,1);
	\draw[fill=gray!50] (2,-2)--(5,-5)--(8,-6)--(3,-1);
	\draw[thick] (2,2) -- (5,5);
	\draw[thick] (2,-2) -- (5,-5);
		\foreach \x in {2,3,4,5,6}
			\draw (\x,\x) -- (\x+1,\x-1);
		\foreach \y in {2,3,4,5,6}
			\draw (\y,-\y) -- (\y+1,-\y+1);
		\foreach \x in {3,4,...,7}
			\draw[dashed,<<-] (\x,\x-2) -- (\x+1,\x-3);
		\foreach \y in {3,4,...,7}
			\draw[dashed,<<-] (\y,-\y+2) -- (\y+1,-\y+3);
	\foreach \x in {1,2}
		\draw[dashed,<<-] (\x,2-\x) -- (\x+1,3-\x);
	\foreach \x in {1,2}
		\draw[dashed,<<-] (\x,-2+\x) -- (\x+1,-3+\x);
	\draw[thick] (0,0) -- (1,1)--(2,0)--(1,-1)--(0,0);
	\draw[thick] (3,1) -- (8,6);
	\draw[thick] (3,-1) -- (8,-6);
	\draw[thick] (4,0) -- (9,5);
	\draw[thick] (4,0) -- (9,-5);
		\foreach \x in {5,6,7,8,9}
			\draw (\x-1,5-\x) -- (\x+4,10-\x);
		\foreach \y in {5,6,7,8,9}
			\draw (\y-1,\y-5) -- (\y+4,\y-10);
\draw[->] (0,0) -- (0,1);
\draw[->] (0,0) -- (1,0);
\draw (0,0) node [left] {\tiny $(0,-\delta)$};
	\draw (2,-2) node [left] {\tiny $(\frac{-\lambda+\delta}{4},\frac{-3\lambda-\delta}{4})$};
	\draw (2,2) node [left] {\tiny $(\frac{-\lambda-\delta}{4},\frac{3\lambda-\delta}{4})$};
	\draw (4,0) node [right] {\tiny $(-\frac{\lambda}{2},\frac{\delta}{2})$};
	
		\draw[right=1pt] (9,0) node {\small $k$};
	\draw[right=1pt] (0,4) node {\small $l$};
	\draw[->] (0,0) -- (0,4);
	\draw[->] (0,0) -- (9,0);
\end{tikzpicture}
\end{figure}

\subsubsection{The Weyl chamber $III_1$}
When the character $\chi_{\delta,\lambda}$ lies in the Weyl chamber $III_1$:
\[
	\lambda+\delta\leq-2, \lambda>0.
\]
Define the submodule $V_2$ of $I(\chi_{\delta, \lambda})$ as a direct sum of the two spaces:
\begin{align*}
	V_{\mathrm{disc}+} (w_{\alpha_1}\chi_{\delta, \lambda})&= U(\frakg)
	\tau^{(0,-\delta)}\\
	V_{\mathbb{H}} (w_{\alpha_1}\chi_{\delta, \lambda}) &= U(\frakg) \tau^{(\frac{\lambda-\delta}{4},\frac{-3\lambda-\delta}{4})}.
\end{align*}
These subspaces form a composition series of $I(\chi_{\delta, \lambda})$:
\[
V_{\mathbb{H}}\oplus V_{\mathrm{disc}+}\xhookrightarrow[V_0/(V_{\mathbb{H}}\oplus V_{\mathrm{disc}+}) = Q_+]{\iota} V_0= I(\chi_{\delta,\lambda})
\]
where
\begin{align*}
V_{\mathbb{H}}(w_{\alpha_1}\chi_{\delta, \lambda})&=\bigoplus_{\substack{(k,l)\in \mathtt{LatticeCond}\\k-l\geq \lambda-\delta}} \tau^{(\frac{k}{2},-\delta+\frac{3l}{2})}\\
Q_{+}(w_{\alpha_1}\chi_{\delta, \lambda})&=\bigoplus_{\substack{(k,l)\in \mathtt{LatticeCond}\\k-l<\lambda-\delta\\k-l\geq-\lambda-\delta}} \tau^{(\frac{k}{2},-\delta+\frac{3l}{2})}\\
V_{\mathrm{disc}+}(w_{\alpha_1}\chi_{\delta, \lambda}) &= \bigoplus_{\substack{(k,l)\in \mathtt{LatticeCond}\\k-l < -\lambda-\delta}} \tau^{(\frac{k}{2}, -\delta+\frac{3l}{2})}.
\end{align*}
For $(\delta,\lambda) = (-6,2)$, the regions representing the $K$-types of these modules are shown in the picture below:
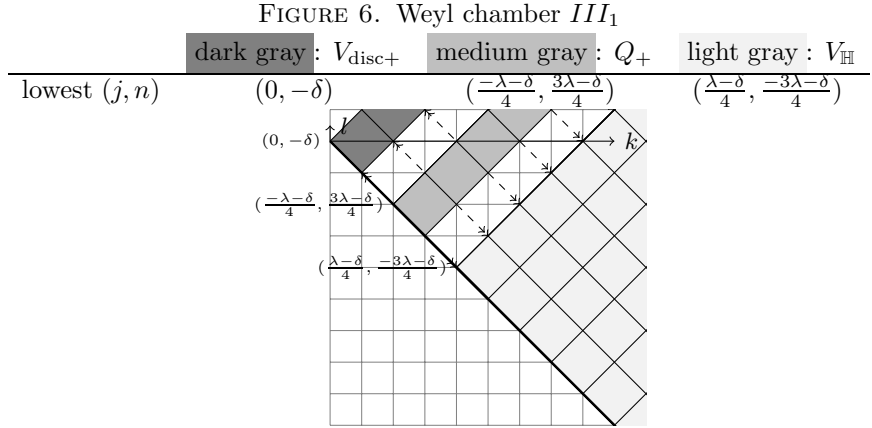
\begin{figure}[h]
\centering
 \caption{Weyl chamber $III_1$}
 \begin{tabular}{cccc}
 &\colorbox{gray}{dark gray}: $V_{\mathrm{disc}+}$&\colorbox{gray!50}{medium gray}: $Q_+$&\colorbox{gray!10}{light gray}: $V_\mathbb{H}$\\\hline
 lowest $(j,n)$ &$(0,-\delta)$&$(\frac{-\lambda-\delta}{4},\frac{3\lambda-\delta}{4})$&$(\frac{\lambda-\delta}{4},\frac{-3\lambda-\delta}{4})$
 \end{tabular}
\begin{tikzpicture}
[scale=0.42]
	\clip (-3.1,-5) rectangle (10,5);
	\draw[step=1,gray,very thin] (0,-5) grid (10,5);
	\draw[fill=gray!10] (4,0) -- (14,-10) -- (14,10);
	\draw[fill=gray!50] (2,2)--(5,5)--(8,6)--(3,1);
	\draw[fill=gray!100] (0,4)--(1,3)--(6,8)--(5,9);
	\foreach \x in {1,2,...,9}
			\draw[dashed,<<-] (\x,\x+2) -- +(1,-1);
		\foreach \y in {3,4,...,9}
			\draw[dashed,->>] (\y,\y-2) -- +(1,-1);
	\foreach \x in {0,1,...,9}
			\draw (\x,\x+4) -- +(1,-1);
	\foreach \x in {2,3,...,9}
			\draw (\x,\x) -- +(1,-1);
	\foreach \x in {4,5,...,10}
			\draw (\x,\x-4) -- +(6,-6);
	\foreach \x in {4,5,...,10}
			\draw (\x,-\x+4) -- +(6,6);
	\draw[thick] (0,4) -- +(5,-5);
	\draw[thick] (1,-1+4) -- +(6,-6);
	\draw[thick] (0,0+4) -- +(5,-5);
	\draw[thick] (2,-2+4) -- +(7,-7);
	\draw[thick] (3,-3+4) -- +(8,-8);
	\draw[thick] (4,-4+4) -- +(9,-9);
	\foreach \x in {4,5,...,8}
		\draw (\x, \x-8+4) -- +(9,9);
	\foreach \x in {5,6,...,8}
		\draw (\x, -\x+4) -- +(9,-9);
	\foreach \x in {0,1,...,8}
			\draw (\x,\x+4) -- +(1,1);
	\foreach \x in {2,3,...,8}
			\draw (\x,\x-4+4) -- +(1,1);
\draw (0,4) node [left] {\tiny $(0,-\delta)$};
	\draw (4,0) node [left] {\tiny $(\frac{\lambda-\delta}{4},\frac{-3\lambda-\delta}{4})$};
	\draw (2,2) node [left] {\tiny $(\frac{-\lambda-\delta}{4},\frac{3\lambda-\delta}{4})$};
	\draw[right=1pt] (9,4) node {\small $k$};
	\draw[right=1pt] (0,4.5) node {\small $l$};
	\draw[->] (0,4) -- (0,4+0.5);
	\draw[->] (0,4) -- (9,4);
    \end{tikzpicture}
\end{figure}

\subsubsection{The Weyl chamber $III_2$}
The character $\chi_{\delta,\lambda}$ lying in Weyl chamber $III_2$ satisfies the inequality: \[\lambda < 0, \lambda-\delta \geq 2.\]There exists a $(\frakg,K)$-submodule submodule $Q_+$ of $I(\chi_{\delta, \lambda})$ defined as:
\begin{align*}
	Q_{+} (w_{\alpha_2}w_{\alpha_1}\chi_{\delta, \lambda})&= U(\frakg)\tau^{(\frac{\lambda-\delta}{4}, \frac{-3\lambda-\delta}{4})}.
%	V_{1} &= \bigoplus_{\substack{(k,l)\in \mathtt{LatticeCond}}} \tau^{(\frac{-\lambda+\delta-2}{4}+\frac{k}{2}, \frac{3\lambda+\delta+6}{4}+\frac{3l}{2})}
\end{align*}
This subspace form a composition series of $I(\chi_{\delta, \lambda})$:
\[
	Q_{+}\xhookrightarrow[V_1/Q_{+} = V_{\mathbb{H}}\oplus V_{\mathrm{disc}+}]{\iota}  V_0= I(\chi_{\delta,\lambda})
\]
where
\begin{align*}	
V_{\mathbb{H}}(w_{\alpha_2}w_{\alpha_1}\chi_{\delta, \lambda})&=\bigoplus_{\substack{(k,l)\in \mathtt{LatticeCond}\\k-l\geq-\lambda-\delta}} \tau^{(\frac{k}{2},-\delta+\frac{3l}{2})}\\
	Q_{+}(w_{\alpha_2}w_{\alpha_1}\chi_{\delta, \lambda}) &= \bigoplus_{\substack{(k,l)\in \mathtt{LatticeCond}\\k-l <-\lambda-\delta\\k-l\geq \lambda-\delta}} \tau^{(\frac{k}{2},-\delta+\frac{3l}{2})}\\
V_{\mathrm{disc}+}(w_{\alpha_2}w_{\alpha_1}\chi_{\delta, \lambda})&=\bigoplus_{\substack{(k,l)\in \mathtt{LatticeCond}\\k-l<\lambda-\delta}} \tau^{(\frac{k}{2},-\delta+\frac{3l}{2})}.
\end{align*}
For $(\delta,\lambda) = (-6,-2)$, the regions representing the $K$-types of these modules are shown in the picture below:
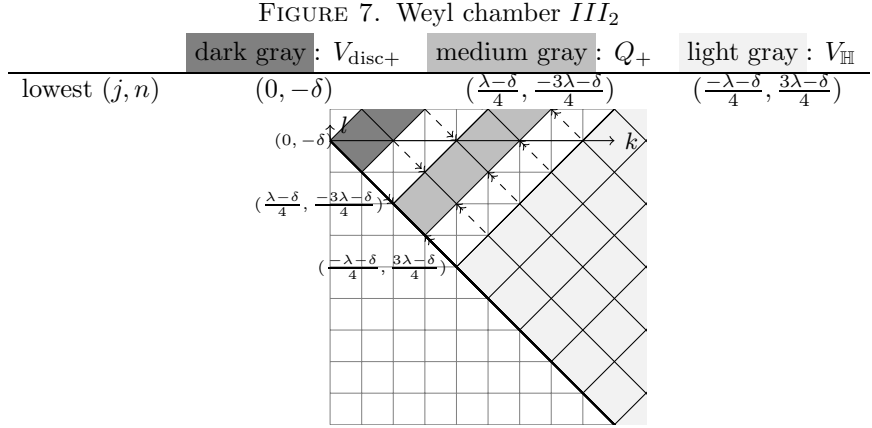
\begin{figure}[h]
\centering
 \caption{Weyl chamber $III_2$}
 \begin{tabular}{cccc}
& \colorbox{gray}{dark gray}: $V_{\mathrm{disc}+}$&\colorbox{gray!50}{medium gray}: $Q_+$&\colorbox{gray!10}{light gray}: $V_\mathbb{H}$\\\hline
lowest $(j,n)$ & $(0,-\delta)$ & $(\frac{\lambda-\delta}{4},\frac{-3\lambda-\delta}{4})$ & $(\frac{-\lambda-\delta}{4},\frac{3\lambda-\delta}{4})$
 \end{tabular}
\begin{tikzpicture}[scale=0.42]
	\clip (-3.1,-5) rectangle (10,5);
	\draw[step=1,gray,very thin] (0,-5) grid (10,5);
	\draw[fill=gray!10] (4,0) -- (14,-10) -- (14,10);
	\draw[fill=gray!50] (2,2)--(5,5)--(8,6)--(3,1);
	\draw[fill=gray!100] (0,4)--(1,3)--(6,8)--(5,9);
	\foreach \x in {1,2,...,9}
			\draw[dashed,->>] (\x,\x+2) -- +(1,-1);
		\foreach \y in {3,4,...,9}
			\draw[dashed,<<-] (\y,\y-2) -- +(1,-1);
	\foreach \x in {0,1,...,9}
			\draw (\x,\x+4) -- +(1,-1);
	\foreach \x in {2,3,...,9}
			\draw (\x,\x) -- +(1,-1);
	\foreach \x in {4,5,...,10}
			\draw (\x,\x-4) -- +(6,-6);
	\foreach \x in {4,5,...,10}
			\draw (\x,-\x+4) -- +(6,6);
	\draw[thick] (0,4) -- +(5,-5);
	\draw[thick] (1,-1+4) -- +(6,-6);
	\draw[thick] (0,0+4) -- +(5,-5);
	\draw[thick] (2,-2+4) -- +(7,-7);
	\draw[thick] (3,-3+4) -- +(8,-8);
	\draw[thick] (4,-4+4) -- +(9,-9);
	\foreach \x in {4,5,...,8}
		\draw (\x, \x-8+4) -- +(9,9);
	\foreach \x in {5,6,...,8}
		\draw (\x, -\x+4) -- +(9,-9);
	\foreach \x in {0,1,...,8}
			\draw (\x,\x+4) -- +(1,1);
	\foreach \x in {2,3,...,8}
			\draw (\x,\x-4+4) -- +(1,1);

\draw (0.4,4) node [left] {\tiny $(0,-\delta)$};
	\draw (4,0) node [left] {\tiny $(\frac{-\lambda-\delta}{4},\frac{3\lambda-\delta}{4})$};
	\draw (2,2) node [left] {\tiny $(\frac{\lambda-\delta}{4},\frac{-3\lambda-\delta}{4})$};
	
	\draw[right=1pt] (9,4) node {\small $k$};
	\draw[right=1pt] (0,4.5) node {\small $l$};
	\draw[->] (0,4) -- (0,4+0.5);
	\draw[->] (0,4) -- (9,4);
    \end{tikzpicture}
\end{figure}

\section{The Intertwining Operator}\label{intertwinesection}
We will prove the Theorem \ref{LongSU21} in this section. The long intertwining operator of the principal series $I(\chi_{\delta,\lambda})$
\[
	A(w_0 ,\chi_{\delta,\lambda}) f(g) = \int_{\overline{N}\cap w^{-1}N w} f(g w_0  \overline{n})\dd \overline{n}
\]
maps each vector $f\in I(\chi_{\delta,\lambda})$ to $A(w_0 ,w_0\chi_{\delta,\lambda}) f\in I(w_0 \chi_{\delta,\lambda})$. We are going to show that this operator acts diagonally on the basis elements $W^{(j,n)}_{m_1,m_2}$ with a closed-form matrix coefficient
\begin{align*}
&\left[A(w,\delta,\lambda)\right]_{m_1} =
	\frac{\pi^2 2^{-\lambda -1}\Gamma (\lambda)}{\Gamma \left(1-\frac{\lambda -\delta }{2}\right)\Gamma \left(1-\frac{\lambda+\delta}{2}\right)}\frac{ \Gamma \left(j+m_1-\frac{\lambda+\delta }{2}+1\right)\Gamma \left(j-m_1-\frac{\lambda-\delta }{2}+1\right)}{
   \Gamma \left(j+m_1+\frac{\lambda-\delta}{2}+1\right) \Gamma \left(j-m_1+\frac{\lambda+\delta }{2}+1\right)}.
\end{align*}
We will start by calculating the Iwasawa decomposition of an element of $w_0\bar{n}$. Since the group $SU(2,1)$ has rank 1, there is only one Weyl group element $w_0 = \diag(-1,-1,1)$ as the reflection of the restricted root system $W(\frakg,\fraka)$. The intersection $\overline{N}\cap w_0 ^{-1}N w_0$ is the set of matrices
\begin{equation}
	\overline{N}\cap w_0 ^{-1}N w_0  =\overline{N}= \{\ppp_{\alpha_1+\alpha_2}\begin{spmatrix}1&0&0\\\sqrt{2}z&1&0\\|z|^2-2\ii w & \sqrt{2}\bar{z} & 1\end{spmatrix}|z\in\CC,w\in \RR\}.
\end{equation}
The matrix $\ppp_{\alpha_1+\alpha_2}\begin{spmatrix}1&0&0\\\sqrt{2}z&1&0\\|z|^2-2\ii w & \sqrt{2}\bar{z} & 1\end{spmatrix}$ has an Iwasawa decomposition in the Lie group $SU(2,1)$:
\begin{align}
	\ppp_{\alpha_1+\alpha_2}&\begin{spmatrix}1&0&0\\\sqrt{2}z&1&0\\|z|^2-2\ii w & \sqrt{2}\bar{z} & 1\end{spmatrix} = \begin{spmatrix}
 -\frac{|z|^2-2 \ii w-1}{\sqrt{(|z|^2+1)^2+4 w^2}} & -\frac{2 \bar{z}}{|z|^2-2 \ii
   w+1} & 0 \\
 \frac{2 z}{\sqrt{(|z|^2+1)^2+4 w^2}} & -\frac{|z|^2+2 \ii w-1}{|z|^2-2 \ii
   w+1} & 0 \\
 0 & 0 & \frac{|z|^2-2 \ii w+1}{\sqrt{(|z|^2+1)^2+4 w^2}} \\
\end{spmatrix} \nonumber\\
&\ppp_{\alpha_1+\alpha_2}\left(\diag\left(\sqrt{(|z|^2+1)^2+4 w^2},1,\frac{1}{\sqrt{(|z|^2+1)^2+4 w^2}}\right)\times \right. \nonumber\\
	&\left.\begin{spmatrix}
 1 & \frac{\sqrt{2} \bar{z}}{\left| z\right| ^2-2 \ii w+1} & \frac{\left| z\right| ^2+2 \ii
   w}{(\left| z\right| ^2+1)^2+4 w^2} \\
 0 & 1 & \frac{\sqrt{2} z}{\left| z\right| ^2+2 \ii w+1} \\
 0 & 0 & 1 \\
\end{spmatrix}\right),\label{IwasawaSU21}
\end{align}\\
where the image of $\ppp_{\alpha_1+\alpha_2}$ on the diagonal matrix lies in $\fraka$, and the image of $\ppp_{\alpha_1+\alpha_2}$ on the upper triangular matrix lies in $N$.
Consider a vector $W^{(j,n)}_{m_1,m_2}$ in $I(\chi_{\delta, \lambda})$. According to the Iwasawa decomposition of an element $\bar{n}\in \bar{N}$ in (\ref{IwasawaSU21}), the right translation of $w_0 \bar{n}$ on this vector can be simplified to
\begin{align}
	W^{(j,n)}_{m_1,m_2}(kw_0\bar{n}) = &\left((|z|^2+1)^2+4 w^2\right)^{-\frac{\lambda+2}{2}}\nonumber\\
	&W^{(j,n)}_{m_1,m_2}\left(kw_0 \begin{spmatrix}
 -\frac{|z|^2-2 \ii w-1}{\sqrt{(|z|^2+1)^2+4 w^2}} & -\frac{2 \bar{z}}{|z|^2-2 i
   w+1} & 0 \\
 \frac{2 z}{\sqrt{(|z|^2+1)^2+4 w^2}} & -\frac{|z|^2+2 \ii w-1}{|z|^2-2 \ii
   w+1} & 0 \\
 0 & 0 & \frac{|z|^2-2 \ii w+1}{\sqrt{(|z|^2+1)^2+4 w^2}} \\
\end{spmatrix}\right).\label{IwasawaActionSU21}
\end{align}
Since $w_0=\diag(-1,-1,1)\in K$, we can absorb $w_0$ by writing
\begin{align*}
	&w_0\begin{spmatrix}
 -\frac{|z|^2-2 \ii w-1}{\sqrt{(|z|^2+1)^2+4 w^2}} & -\frac{2 \bar{z}}{|z|^2-2 i
   w+1} & 0 \\
 \frac{2 z}{\sqrt{(|z|^2+1)^2+4 w^2}} & -\frac{|z|^2+2 \ii w-1}{|z|^2-2 \ii
   w+1} & 0 \\
 0 & 0 & \frac{|z|^2-2 \ii w+1}{\sqrt{(|z|^2+1)^2+4 w^2}} \\
\end{spmatrix} =\\
& \begin{spmatrix}
 \frac{|z|^2-2 \ii w-1}{\sqrt{(|z|^2+1)^2+4 w^2}} & \frac{2 \bar{z}}{|z|^2-2 i
   w+1} & 0 \\
 -\frac{2 z}{\sqrt{(|z|^2+1)^2+4 w^2}} & \frac{|z|^2+2 \ii w-1}{|z|^2-2 \ii
   w+1} & 0 \\
 0 & 0 & \frac{|z|^2-2 \ii w+1}{\sqrt{(|z|^2+1)^2+4 w^2}} \\
\end{spmatrix}.
\end{align*} If we set $z = x+\ii y$, the Haar measure on $N$ is given by $\dd x\dd y\dd w$. Therefore, to understand the intertwining operator $A(w,\delta,\lambda)$, it suffices to compute the singular integral
\begin{align}
	&\int_{\CC\times\RR}\left((|z|^2+1)^2+4 w^2\right)^{-\frac{\lambda+2}{2}}\times\nonumber\\
	&W^{(j,n)}_{m_1,m_2}\begin{spmatrix}
 \frac{|z|^2-2 \ii w-1}{\sqrt{(|z|^2+1)^2+4 w^2}} & \frac{2 \bar{z}}{|z|^2-2 i
   w+1} & 0 \\
 -\frac{2 z}{\sqrt{(|z|^2+1)^2+4 w^2}} & \frac{|z|^2+2 \ii w-1}{|z|^2-2 \ii
   w+1} & 0 \\
 0 & 0 & \frac{|z|^2-2 \ii w+1}{\sqrt{(|z|^2+1)^2+4 w^2}} \\
\end{spmatrix} \dd x \dd y\dd w\label{FormOfSU21Integral}
\end{align}
for every $K$-type $(j,n)$ and all indices in $-j\leq m_1,m_2\leq j$. We will calculate the integral in the domain of $\lambda$ where it converges, and deduce the validity of the formula in Theorem 1.3 by analytic continuation.  According to the following theorem from \cite{VoganWallachIntertwining}, the long intertwining operator $A(w_0,\chi_{\delta,\lambda})$ depends meromorphically on $\lambda\in\frakh_\CC^*$:
\begin{theorem}
	There exist polynomial maps $b_\delta:\fraka_\CC^*\longrightarrow\CC$ and $D_\delta:\fraka_\CC^*\longrightarrow U(\frakg_\CC)^K$, such that for $f\in I_P(\chi_{\delta,\lambda})$, if $\lambda$ satisfies $\mathrm{Re}\langle\lambda,\alpha_i\rangle\geq c_{\delta}$ for $\alpha\in\Sigma^+(\frakg,\fraka)$, we have
	\[
		b_{\delta}(\lambda)A(P|\bar{P},\chi_{\delta,\lambda}) f= A(P|\bar{P},\chi_{\delta,\lambda+4\rho})\pi_{P}(\chi_{\delta,\lambda+4\rho})(D_{\delta})f.
	\]
\end{theorem}
According to this theorem, there exists a number $c_\delta>0$, such that the integral (\ref{FormOfSU21Integral}) converges if $\mathrm{Re}\lambda>c_\delta$. Based on the definition of Wigner $D$-functions (3.13), the integrand can be expressed as a hypergeometric sum
\begin{align}
	&\left((|z|^2+1)^2+4 w^2\right)^{-\frac{\lambda+2}{2}}W^{(j,n)}_{m_1,m_2}\begin{spmatrix}
 \frac{|z|^2-2 \ii w-1}{\sqrt{(|z|^2+1)^2+4 w^2}} & \frac{2 \bar{z}}{|z|^2-2 i
   w+1} & 0 \\
 -\frac{2 z}{\sqrt{(|z|^2+1)^2+4 w^2}} & \frac{|z|^2+2 \ii w-1}{|z|^2-2 \ii
   w+1} & 0 \\
 0 & 0 & \frac{|z|^2-2 \ii w+1}{\sqrt{(|z|^2+1)^2+4 w^2}} \\
\end{spmatrix}\nonumber\\
=&c^j_{m_1}c^j_{m_2}\sum_{p=\max(0,m_1-m_2)}^{\min(j-m_2,j+m_1)}\frac{(-1)^{p}2^{-m_1+m_2+2p}}{(j+m_1-p)!p!(m_2-m_1+p)!(j-m_2-p)!}\times\nonumber\\&\omega^{(j,n)}_{m_1,m_2}(p;z,w)\label{HyperSumSU21}
\end{align}
where the function $\omega^{(j,n)}_{m_1,m_2}(p;z,w)$ is defined as a function in $z\in\CC, w\in\RR$
\begin{align}
	&\omega^{(j,n)}_{m_1,m_2}(p;z,w) = z^p\bar{z}^{-m_1+m_2+p}(-1+|z|^2+2\ii w)^{j+m_1-p}(-1+|z|^2-2\ii w)^{j-m_2-p}\nonumber\\
	&(1+|z|^2-2\ii w)^{\frac{-2j-m_2+n-\lambda-2}{2}}(1+|z|^2+2\ii w)^{\frac{-2j+m_2-n-\lambda-2}{2}}\label{FunctionW}
\end{align}
which can be factored into polynomial functions in $z$ and $w$. Noting that since $1+|z|^2>0$, the complex number $1+|z|^2\pm2\ii w$ lies in the right half plane, we can always take a branch cut of the power functions in (\ref{FunctionW}) such that the value of $(1+|z|^2\mp 2\ii w)^{\frac{-2j\mp m_2\pm n-\lambda-2}{2}}$ when $z=0, w=0$ is 1.\\

In order to compute the integral (\ref{FormOfSU21Integral}), it suffices to integrate on each summand $\omega^{(j,n)}_{m_1,m_2}(p;z,w)$ over $\CC\times\RR$. We can change the rectangular coordinate $z = x+\ii y$ to the polar coordinate $z=r e^{\ii \theta}$, and by (\ref{FunctionW}), we have
\begin{align}
	&\int_{\CC\times\RR}\omega^{(j,n)}_{m_1,m_2}(p;z,w)\dd x\dd y\dd w\nonumber\\
	&=\int_{0}^{\infty} r^{-m_1+m_2+2p+1}\dd r\int_{-\infty}^{\infty}(-1+r^2+2\ii w)^{j+m_1-p}(-1+r^2-2\ii w)^{j-m_2-p}\nonumber\\
	&(1+r^2-2\ii w)^{\frac{-2j-m_2+n-\lambda-2}{2}}(1+r^2+2\ii w)^{\frac{-2j+m_2-n-\lambda-2}{2}}\dd w\int_{0}^{2\pi}e^{\ii(m_1-m_2)\theta}\dd \theta\nonumber\\
	&=2\pi\delta_{m_1,m_2}\int_{0}^{\infty} r^{-m_1+m_2+2p+1}\left(\int_{-\infty}^{\infty}(-1+r^2+2\ii w)^{j+m_1-p}(-1+r^2-2\ii w)^{j-m_2-p}\right.\nonumber\\
	&\left.(1+r^2-2\ii w)^{\frac{-2j-m_2+n-\lambda-2}{2}}(1+r^2+2\ii w)^{\frac{-2j+m_2-n-\lambda-2}{2}}\dd w\right)\dd r\label{PolarIntegral}
\end{align}
From the last line of the calculation above, we notice that the matrix \[(\int_{\CC\times\RR}\omega^{(j,n)}_{m_1,m_2}(p;z,w)\dd x\dd y\dd w)_{-j\leq m_1,m_2\leq j}\] is diagonal due to the appearance of $\delta_{m_1,m_2}$, i.e. its entries are nonzero if and only if $m_1=m_2$. Therefore, the intertwining operator $A(w,\delta,\lambda)$ acts diagonally on each $K$-type, and write the diagonal entries as $\left[A(w,\delta,\lambda)\right]_{m_1} = \langle W^{(j,n)}_{m_1,*}, A(w,\delta,\lambda)W^{(j,n)}_{m_1,*}\rangle$. We define $\tilde{\omega}^{(j,n)}_{m_1}(r,w)$ as the inner integrand of the integral above:
\begin{align*}
	\tilde{\omega}^{(j,n)}_{m_1}(r,w)=&(-1+r^2+2\ii w)^{j+m_1-p}(-1+r^2-2\ii w)^{j-m_1-p}
	(1+r^2-2\ii w)^{\frac{-2j-m_1+n-\lambda-2}{2}}\\
	&(1+r^2+2\ii w)^{\frac{-2j+m_1-n-\lambda-2}{2}},
\end{align*}
so that
\[
\int_{\CC\times\RR}\omega^{(j,n)}_{m_1,m_2}(p;z,w)\dd x\dd y\dd w = 2\pi \delta_{m_1,m_2}\int_{0}^{\infty} r^{-m_1+m_2+2p+1}\left(\int_{-\infty}^{\infty}\tilde{\omega}^{(j,n)}_{m_1}(r,w)\dd w\right)\dd r .
\]
Using the notation which we have just introduced, the diagonal elements of the intertwining operator can be expressed as
\begin{align}
	\left[A(w,\delta,\lambda)\right]_{m_1} &= 2\pi \sum_{p=0}^{\min(j-m_1,j+m_1)}\begin{spmatrix}
	j+m_1\\p
	\end{spmatrix}\begin{spmatrix}
	j-m_1\\p
	\end{spmatrix}(-4)^p\times\nonumber\\ &\int_0^\infty r^{2p+1}\left(\int_{-\infty}^{\infty}\tilde{\omega}^{(j,n)}_{m_1}(r,w)\dd w\right)\dd r \label{SU21Integral:00}
\end{align}
Apply the change of variables  \[\begin{matrix}m_1=\frac{n+\delta}{3},&j=\frac{k}{2},&n=-\delta+\frac{3}{2}l\end{matrix}\] as we used when describing the structure of $SU(2,1)$ principal series in Section {\ref{decomp}}, the integrand $\tilde{\omega}^{(j,n)}_{m_1}(r,w)$ inside can be replaced by another function labelled by $k,l$
\begin{align}
	\tilde{\omega}_{k,l}(r,w)=&(-1+r^2+2\ii w)^{\frac{k+l}{2}-p}(-1+r^2-2\ii w)^{\frac{k-l}{2}-p}(1+r^2-2\ii w)^{-\frac{k-l+\lambda+\delta+2}{2}}\nonumber\\
    &(1+r^2+2\ii w)^{-\frac{k+l+\lambda-\delta+2}{2}}\label{SU21Integral:0}.
\end{align}
The integral $\int_{-\infty}^{\infty}\tilde{\omega_{k,l}}(r,w)\dd w$ converges for $\lambda>0$. In order to calculate the integral, we need to reorganize the factors for  $\tilde{\omega}_{k,l}(r,w)$ to a simpler form. By applying the change of variable from $w$ to $\frac{1}{2}(1+r^2) w$, the integral $\int_{-\infty}^{\infty}\tilde{\omega_{k,l}}(r,w)\dd w$ can be rewritten as follows:
\begin{align}
	&\int_{-\infty}^{\infty}(-1+r^2+2\ii w)^{\frac{k+l}{2}-p}(-1+r^2-2\ii w)^{\frac{k-l}{2}-p}(1+r^2-2\ii w)^{-\frac{k-l+\lambda+\delta+2}{2}}\nonumber\\
	&(1+r^2+2\ii w)^{-\frac{k+l+\lambda-\delta+2}{2}}\dd w\nonumber\\
	=& \frac{1}{2} (1+r^2)^{-1-2p-\lambda}\int_{-\infty}^{\infty}\left(\frac{-1+r^2}{1+r^2}+\ii w\right)^{\frac{k+l}{2}-p}\left(\frac{-1+r^2}{1+r^2}-\ii w\right)^{\frac{k-l}{2}-p}\nonumber\\
	&(1-\ii w)^{-\frac{k-l+\lambda+\delta+2}{2}}(1+\ii w)^{-\frac{k+l+\lambda-\delta+2}{2}}\dd w.\label{SU21Integral:1}
\end{align}
Since in the original summation of Wigner $D$-functions, when $m_1=m_2=m$, $0\leq p\leq \min(j-m,j+m)=\min(\frac{k-l}{2},\frac{k+l}{2})$, the exponents of the first two factors \[\left(\frac{-1+r^2}{1+r^2}+\ii w\right)^{\frac{k+l}{2}-p}\left(\frac{-1+r^2}{1+r^2}-\ii w\right)^{\frac{k-l}{2}-p}\] in the integrand (\ref{SU21Integral:1}) are non-negative integers, we can thus reorganize the terms inside of the parenthesis and expand using binomial theorem:
\begin{align}
	&\left(\frac{-1+r^2}{1+r^2}+\ii w\right)^{\frac{k+l}{2}-p}\left(\frac{-1+r^2}{1+r^2}-\ii w\right)^{\frac{k-l}{2}-p}\nonumber\\
	=&\left(\left(\frac{-1+r^2}{1+r^2}+1\right)-1+\ii w\right)^{\frac{k+l}{2}-p}\left(\left(\frac{-1+r^2}{1+r^2}-1\right)+1-\ii w\right)^{\frac{k-l}{2}-p}\nonumber\\
	=& \sum_{K_1,K_2}(-1)^{\frac{k+l}{2}-p-K2}\begin{spmatrix}\frac{k-l}{2}-p\\K_1\end{spmatrix}\begin{spmatrix}\frac{k+l}{2}-p\\K_2\end{spmatrix}\left(\frac{-1+r^2}{1+r^2}-1\right)^{K_1}\left(\frac{-1+r^2}{1+r^2}+1\right)^{K_2}\nonumber\\
	&(1-\ii w)^{k-2p-K_1-K_2}.\label{SU21Integral:2}
\end{align}
Combining the factor $(1-\ii w)^{k-2p-K_1-K_2}$ with (\ref{SU21Integral:1}), the integral in (\ref{SU21Integral:00}) which depends on $w$ becomes
\begin{align}
	&\int_{-\infty}^{\infty}(1-\ii w)^{-\frac{k-l+(\lambda+\delta+2)}{2}+k-2p-K_1-K_2}(1+\ii w)^{-\frac{k+l+(\lambda-\delta+2)}{2}}\dd w\nonumber\\
	=&\frac{2^{-K_1-K_2-2p-\lambda}\pi \Gamma(1+K_1+K_2+2p+\lambda)}{\Gamma\left(\frac{k+l+\lambda-\delta}{2}+1\right)\Gamma\left(-\frac{k+l-\lambda-\delta}{2}+1+K_1+K_2+2p\right)}.
\end{align}
Then we can take care of the integral which depends on $r$ in (\ref{SU21Integral:00}):
\begin{align}
	&\int_{0}^\infty (1+r^2)^{-1-2p-\lambda}\left(\frac{-1+r^2}{1+r^2}-1\right)^{K_1}\left(\frac{-1+r^2}{1+r^2}+1\right)^{K_2}r^{2p+1}\dd r\nonumber\\
	=&(-1)^{K_1}2^{K_1+K_2-1}\frac{\Gamma(1+K_2+p)\Gamma(K_1+p+\lambda)}{\Gamma(1+K_1+K_2+2p+\lambda)}\label{SU21Integral:3}
\end{align}
Putting (\ref{SU21Integral:0})-(\ref{SU21Integral:1}) back into the intertwining operator $\left[A(w,\delta,\lambda)\right]_{m_1}$ integral (\ref{SU21Integral:00}) and applying the change of indices in $j,m_1$ to $k,l$, the summation in (\ref{SU21Integral:00}) becomes a sum over $\Gamma$-functions and binomial coefficients. We can utilize a trick by changing all the binomial coefficients into their $\Gamma$ function expressions, and group the $\Gamma$-factors in the following way:
\begin{align*}
	&\left[A(w,\delta,\lambda)\right]_{m_1} =  2^{-\lambda-1}(-1)^{\frac{k+l}{2}}\pi^2 \frac{\Gamma \left(\frac{k+l+2}{2}\right) \Gamma \left(\frac{
   k-l+2}{2}\right)}{ \Gamma \left(\frac{k+l+\lambda-\delta}{2}+1\right)}\\
   &\sum_{\substack{K_1,K_2\geq 0\\p\geq 0}}
   \frac{(-1)^{K_1+K_2}}{\Gamma \left(K_1+1\right) \Gamma
   \left(K_2+1\right) \Gamma \left(K_1+K_2+2p-\frac{k+l-\lambda-\delta}{2}+1\right)}\\
   &\frac{\Gamma
   \left(p+\lambda +K_1\right)\Gamma \left(1+K_2+p\right) }{\Gamma (p+1)^2 \Gamma
   \left(\frac{k-l}{2}-p-K_1+1\right)\Gamma
   \left(\frac{k+l}{2}-p-K_2+1\right) }
\end{align*}
Reorganizing the $\Gamma$-functions into multinomial coefficients
\[
	\begin{spmatrix}n\\n_1,\ldots,n_r\end{spmatrix} = \frac{n!}{n_1!\cdots n_r!}
\]and adding auxiliary $\Gamma$ factors as required, we get:
\begin{align}
	&\left[A(w,\delta,\lambda)\right]_{m_1} =  2^{-\lambda-1}(-1)^{\frac{k+l}{2}}\pi^2 \frac{\Gamma \left(\frac{k+l+2}{2}\right) \Gamma \left(\frac{k-l+2}{2}\right)\Gamma(\lambda)}{ \Gamma \left(\frac{k+l-\delta +\lambda }{2}+1\right)\Gamma \left(\frac{k-l+\delta +\lambda}{2}+1\right)}\nonumber\\&\sum_{\substack{K_1,K_2\geq 0\\p\geq 0}}(-1)^{K_1+K_2}
   \begin{spmatrix}K_1+p+\lambda-1\\\lambda-1,K_1,p\end{spmatrix}
   \begin{spmatrix}K_2+p\\p\end{spmatrix}\begin{spmatrix}\frac{k-l+\lambda+\delta}{2}\\\frac{k-l}{2}-p-K_1,\frac{k+l}{2}-p-K_2,-\frac{k+l-\lambda-\delta}{2}+K_1+K_2+2p\end{spmatrix}.\label{SU21Binomial}
\end{align}
Noticing that the $y^{K_1+p}z^{K_2+p}$-th multinomial coefficient of the following function in $y,z$:
\[
	y^{\frac{k-l}{2}}z^{\frac{k+l}{2}}\left(1-\frac{1}{y}-\frac{1}{z}\right)^{\frac{k-l+\lambda+\delta}{2}}
\]
is the multinomial coefficient \[(-1)^{k+K_1+K_2}\begin{spmatrix}\frac{k-l+\lambda+\delta}{2}\\\frac{k-l}{2}-p-K_1,\frac{k+l}{2}-p-K_2,-\frac{k+l-\lambda-\delta}{2}+K_1+K_2+2p\end{spmatrix}\]in the third factor of each summand, we can further take $y = (1+s)(1+t)$ and $z =1+1/t$. We consider the multinomial expansion of the function
\begin{align}
	&(1+s)^{\frac{k-l}{2}+\lambda-1}(1+t)^{\frac{k-l}{2}}(1+1/t)^{\frac{k+l}{2}}\left(1-\frac{1}{(1+s)(1+t)}-\frac{1}{(1+1/t)}\right)^{\frac{k-l+\lambda+\delta}{2}}\nonumber\\
	=&\sum_{\substack{\kappa_1,\kappa_2\in\ZZ}}(-1)^{k+\kappa_1+\kappa_2}\begin{spmatrix}\frac{k-l+\lambda+\delta}{2}\\\frac{k-l}{2}-\kappa_1,\frac{k+l}{2}-\kappa_2,-\frac{k+l-\lambda-\delta}{2}+\kappa_1+\kappa_2\end{spmatrix}(1+s)^{\kappa_1+\lambda-1}\nonumber\\&(1+t)^{\kappa_1}(1+1/t)^{\kappa_2},
\end{align}
its coefficient of the term $s^{\lambda-1}t^0$ is given by
\begin{align}
	&\sum_{\substack{\kappa_1,\kappa_2\in\ZZ}}(-1)^{k+\kappa_1+\kappa_2}\begin{spmatrix}\frac{k-l+\lambda+\delta}{2}\\\frac{k-l}{2}-\kappa_1,\frac{k+l}{2}-\kappa_2,-\frac{k+l-\lambda-\delta}{2}+\kappa_1+\kappa_2\end{spmatrix}\begin{spmatrix}\kappa_1+\lambda-1\\\kappa_1\end{spmatrix}\nonumber\\&\sum_{p\in\ZZ} \begin{spmatrix}\kappa_1\\p\end{spmatrix}\begin{spmatrix}\kappa_2\\p\end{spmatrix}
\end{align}
which is $(-1)^k$ times the sum in (\ref{SU21Binomial}) if one changes $K_i+p$ to $\kappa_i$. The sums in  (\ref{SU21Binomial}) are over non-negative integers, which is guaranteed by the non-vanishing of binomial coefficients. Putting these binomial coefficients all together, it is clear that the sum in (\ref{SU21Binomial}) is the constant term coefficient of the function
\begin{equation}
	(-1)^{-k}s^{1+\frac{1}{2} (k-l-\lambda+\delta )} (1+s)^{\frac{\lambda -\delta }{2}-1}t^{-\frac{k+l}{2}} (1+t)^{\frac{k+l-\delta -\lambda}{2}} 
\end{equation}
which by the binomial theorem is equal to
\begin{equation}
	(-1)^{-k}\frac{ \Gamma \left(1+\frac{k+l}{2}-\frac{\lambda +\delta
   }{2}\right) \Gamma \left(\frac{\lambda -\delta }{2}\right)}{\Gamma
   \left(\frac{k+l}{2}+1\right) \Gamma \left(\frac{k-l}{2}+1\right) \Gamma
   \left(1-\frac{\lambda +\delta }{2}\right) \Gamma \left(-\frac{1}{2}
   (k-l)+\frac{\lambda -\delta }{2}\right)}.
\end{equation}
The matrix entries $\left[A(w,\delta,\lambda)\right]_{m_1}$ for the intertwining operator simplifies to the function
\begin{align}
	&\left[A(w,\delta,\lambda)\right]_{m_1} =\nonumber\\
	& \frac{(-1)^{\frac{k-l}{2}}\pi^2}{2^{\lambda+1}}\frac{\Gamma(\lambda)\Gamma \left(\frac{\lambda -\delta }{2}\right)}{\Gamma
   \left(1-\frac{\lambda +\delta }{2}\right) }\frac{\Gamma \left(\frac{k+l-\lambda-\delta+2}{2}\right) }{\Gamma \left(-\frac{k-l-\lambda+\delta}{2}
   \right) \Gamma \left(\frac{k+l +\lambda -\delta+2}{2}\right)\Gamma \left(\frac{k-l+\lambda +\delta +2}{2}\right)}.
\end{align}
If we change the indices $(k,l)$ back to $j, m_1$, we have:
\begin{align*}
&\left[A(w,\delta,\lambda)\right]_{m_1} =
	\frac{\pi ^2 2^{-\lambda -1} (-1)^{
   j-m_1} \Gamma (\lambda ) \Gamma \left(\frac{\lambda -\delta }{2}\right)}{\Gamma \left(1-\frac{\lambda+\delta}{2}\right)}\\&\frac{ \Gamma \left(j+m_1-\frac{\lambda+\delta }{2}+1\right)}{ \Gamma \left(-j+m_1+\frac{\lambda-\delta}{2}\right)
   \Gamma \left(j+m_1+\frac{\lambda-\delta}{2}+1\right) \Gamma \left(j-m_1+\frac{\lambda+\delta }{2}+1\right)}
\end{align*}
We can apply the formula
\[
	\Gamma(z)\Gamma(1-z) = \frac{\pi}{\sin(\pi z)}
\]
we can move one of the $\Gamma$-factors from the denominator to the numerator and vice versa, which gives the final formula for the long intertwining operator entries:
\begin{align}
&\left[A(w,\delta,\lambda)\right]_{m_1} =
	\frac{\pi^2 2^{-\lambda -1}\Gamma (\lambda)}{\Gamma \left(1-\frac{\lambda -\delta }{2}\right)\Gamma \left(1-\frac{\lambda+\delta}{2}\right)}\frac{ \Gamma \left(j+m_1-\frac{\lambda+\delta }{2}+1\right)\Gamma \left(j-m_1-\frac{\lambda-\delta }{2}+1\right)}{
   \Gamma \left(j+m_1+\frac{\lambda-\delta}{2}+1\right) \Gamma \left(j-m_1+\frac{\lambda+\delta }{2}+1\right)}
\end{align}
Thus we have proven Theorem \ref{LongSU21}.

\bibliographystyle{alpha}
\bibliography{su21}

\end{document}